\documentclass[12pt]{article}

\usepackage{amssymb,amsthm,amscd,array}

\let\ssection=\section
\renewcommand{\section}{\setcounter{equation}{0}\ssection}

\newcommand{\bbR}{\mathbb{R}}

\newcommand{\bbC}{\mathbb{C}}

\newcommand{\cC}{{\mathcal{C}}}
\newcommand{\diag}{\mathrm{diag}}
\newcommand{\cD}{{\mathcal{D}}}

\newcommand{\cE}{{\mathcal{E}}}
\newcommand{\Diff}{\mathrm{Diff}}
\newcommand{\Div}{\mathrm{Div}}
\newcommand{\rD}{\mathrm{D}}
\newcommand{\rE}{\mathrm{E}}
\newcommand{\Euler}{{\mathcal E}}
\newcommand{\End}{\mathrm{End}}
\newcommand{\cF}{{\mathcal{F}}}
\newcommand{\gl}{{\mathrm{gl}}}

\newcommand{\rg}{\mathrm{g}}
\newcommand{\rG}{\mathrm{G}}

\newcommand{\rh}{\mathrm{h}}
\newcommand{\cI}{{\mathcal{I}}}
\newcommand{\Id}{\mathrm{Id}}
\newcommand{\Or}{\mathrm{O}}
\newcommand{\Pol}{\mathrm{Pol}}
\newcommand{\cQ}{{\mathcal{Q}}}
\newcommand{\rL}{\mathrm{L}}
\newcommand{\rR}{\mathrm{R}}
\newcommand{\rT}{\mathrm{T}}

\newcommand{\Ric}{\mathrm{Ric}}
\newcommand{\cS}{{\mathcal{S}}}

\newcommand{\SL}{\mathrm{SL}}
\newcommand{\Sl}{\mathrm{sl}}

\newcommand{\se}{\mathrm{e}}
\newcommand{\so}{\mathrm{o}}

\newcommand{\ssp}{\mathrm{sp}}

\newcommand{\Supp}{\mathrm{Supp}}

\newcommand{\Vect}{\mathrm{Vect}}
\newcommand{\vol}{\mathrm{vol}}

\newcommand{\half}{\frac{1}{2}}

\def\mod#1{\left|{#1}\right|}

\begin{document}

%\baselineskip=18pt

%\textwidth=16truecm
%\textheight=24truecm
%\hoffset=-1.5truecm
%\voffset=-2.5truecm

\def\a{\alpha}
\def\b{\beta}
\def\d{\delta}
\def\g{\gamma}
\def\om{\omega}
\def\r{\rho}
\def\s{\sigma}
\def\vfi{\varphi}
\def\l{\lambda}
\def\m{\mu}
\def\implies{\Rightarrow}

\oddsidemargin .1truein
\newtheorem{thm}{Theorem}[section]
\newtheorem{lem}[thm]{Lemma}
\newtheorem{cor}[thm]{Corollary}
\newtheorem{pro}[thm]{Proposition}
\newtheorem{ex}[thm]{Example}
\newtheorem{rmk}[thm]{Remark}
\newtheorem{defi}[thm]{Definition}
%\newremark{ex}[thm]{Example}

%\newtheorem{thm}{Theorem}
%\newtheorem{lem}{Lemma}
%\newtheorem{cor}{Corollary}
%\newtheorem{prop}[thm]{Proposition}
%\newtheorem{definition}{Definition}

\title{Conformally equivariant quantization}

\author{C.~Duval\footnote{mailto:duval@cpt.univ-mrs.fr}\\
{\small Universit\'e de la M\'editerran\'ee and CPT-CNRS}
\and
V.~Ovsienko\footnote{mailto:ovsienko@cpt.univ-mrs.fr}\\
{\small CNRS, Centre de Physique Th\'eorique}
\thanks{CPT-CNRS, Luminy Case 907,
F--13288 Marseille, Cedex 9, FRANCE.
}
}

\date{}

\maketitle

\thispagestyle{empty}

\begin{abstract}
Let $(M,\rg)$ be a pseudo-Riemannian manifold and $\cF_\l(M)$ the space of densities
of degree $\l$ on $M$. We study the space $\cD^2_{\l,\m}(M)$ of second-order
differential operators  from $\cF_\l(M)$ to $\cF_\m(M)$. If $(M,\rg)$ is conformally
flat with signature $p-q$, then $\cD^2_{\l,\m}(M)$ is  viewed as a module over the
group of conformal transformations of $M$. We prove that, for almost all values of
$\m-\l$, the $\Or(p+1,q+1)$-modules $\cD^2_{\l,\m}(M)$ and the space of symbols
(i.e., of second-order polynomials on $T^*M$)  are canonically isomorphic. This
yields a conformally equivariant quantization for quadratic Hamiltonians. We
furthermore show that this quantization map extends to arbitrary pseudo-Riemannian
manifolds and depends only on the conformal class $[\rg]$ of the metric.   As an
example, the quantization of the geodesic flow yields a novel conformally
equivariant Laplace operator on half-densities, as well as the well-known Yamabe
Laplacian. We also recover in this framework the multi-dimensional Schwarzian
derivative of conformal transformations.
\end{abstract}


%\newpage

%%%%%%%%%%%%%%%%%%%%%%%%%%%%%%%%%%%%%%%%%%%%%%%%%%%%%%%%%%%%%%%%%%%%%%%%%%%%%%%%%%%%
%%%%%%%%%%%%%%%%%%%%%%%%%%%%%%%%%%%%%%%%%%%%%%%%%%%%%%%%%%%%%%%%%%%%%%%%%%%%%%%%%%%%
\section{Introduction}
%%%%%%%%%%%%%%%%%%%%%%%%%%%%%%%%%%%%%%%%%%%%%%%%%%%%%%%%%%%%%%%%%%%%%%%%%%%%%%%%%%%%
%%%%%%%%%%%%%%%%%%%%%%%%%%%%%%%%%%%%%%%%%%%%%%%%%%%%%%%%%%%%%%%%%%%%%%%%%%%%%%%%%%%%

The aim of this article is to investigate the relationship between
differential operators on a smooth pseudo-Riemannian manifold $(M,\rg)$ of
signature $p-q$, and the  polynomial functions on its cotangent bundle $T^*M$. 

%%%%%%%%%%%%%%%%%%%%%%%%%%%%%%%%%%%%%%%%%%%%%%%%%%%%%%%%%%%%%%%%%%%%%%%%%%%%%%%%%%%%
%\subsection{Statement of the problem}
%%%%%%%%%%%%%%%%%%%%%%%%%%%%%%%%%%%%%%%%%%%%%%%%%%%%%%%%%%%%%%%%%%%%%%%%%%%%%%%%%%%%

We will consider the space $\cD(M)$ of differential operators on
$C^\infty$-function of~$M$ viewed as a module for the group $\Diff(M)$ of all
diffeomorphisms of $M$. We are, in fact, interested in a two-parameter family of 
modules which can be understood as follows.
Considering that the arguments of these differential operators are, indeed,
tensor densities of, say, weight $\l$, and their values tensor densities of
weight $\m$, we will, hence, deal with a new $\Diff(M)$-module structure
denoted by $\cD_{\l,\m}(M)$.

The natural $\Diff(M)$-module of symbols associated with $\cD_{\l,\m}(M)$ is the
space of fiberwise polynomials on $T^*M$ with values in the $(\m-\l)$-densities over
$M$. Therefore, we have a one-parameter family of $\Diff(M)$-modules, $\cS_\d(M)$,
where $\d=\m-\l$.

%\begin{rmk}{\rm
The modules $\cD_{\l,\m}(M)$ have already been considered in the classic literature
on differential operators and, more recently, in a series of papers
\cite{DO,LMT,GO,LO,Gar,Mat}. The general problem of classification of these
$\Diff(M)$-modules has been solved in these articles. 
%}\end{rmk}

We will be considering the modules of second-order
operators, $\cD_{\l,\m}^2(M)$, and symbols, $\cS_\d^2(M)$.

The main purpose of this article is to define a canonical isomorphism
\begin{equation}
\cQ_{\l,\m}:\cS^2_{\m-\l}(M)\stackrel{\cong}\longrightarrow\cD^2_{\l,\m}(M)
\label{confEquivIsom}
\end{equation}
that satisfies the following properties:

\begin{enumerate}
\item
it is conformally \textit{invariant}, i.e., it depends only on the
conformal class $[\rg]$ of the metric;

\item
in the conformally flat case, it is \textit{equivariant} with respect to
$\Or(p+1,q+1)$, the group of conformal diffeomorphisms.
\end{enumerate}

We will show that the isomorphism (\ref{confEquivIsom}) exists for generic 
$\l$ and $\m$; in the most interesting case
$\l=\m=\half$, it provides a natural quantization of the cotangent bundle
of a pseudo-Riemannian manifold.

In the conformally flat case, the isomorphism (\ref{confEquivIsom}) is 
characterized by the equivariance property and is essentially unique (up to a natural
normalization).

%%%%%%%%%%%%%%%%%%%%%%%%%%%%%%%%%%%%%%%%%%%%%%%%%%%%%%%%%%%%%%%%%%%%%%%%%%%%%%%%%%%%
%\subsection{Existence theorems}
%%%%%%%%%%%%%%%%%%%%%%%%%%%%%%%%%%%%%%%%%%%%%%%%%%%%%%%%%%%%%%%%%%%%%%%%%%%%%%%%%%%%

It worth noticing that for any $\l$ and $\m$ there is no isomorphism
(\ref{confEquivIsom})
equivariant with respect to the full group $\Diff(M)$.

Let us assume that the manifold $M$ is endowed with a flat conformal
structure which enables us to look for a conformally-equivariant quantization
with respect to the group $\Or(p+1,q+1)$ (or its Lie algebra $\so(p+1,q+1)$) where
$\dim(M)=p+q$ acting (locally) on~$M$. 

\begin{thm}\label{confEquivIsomTh}
Given a conformally flat pseudo-Riemannian manifold~$M$ of dimension
$n=p+q\geq2$, 
\hfill\break
(i) there exists an isomorphism (\ref{confEquivIsom}) of $\so(p+1,q+1)$-modules
provided
\begin{equation}
\m-\l\not\in\left\{\frac{2}{n},
\frac{n+2}{2n},
1,
\frac{n+1}{n},
\frac{n+2}{n}\right\}.
\label{Reson}
\end{equation}
\hfill\break
(ii) For every $\l$ and $\m$ as in (\ref{Reson}), this isomorphism is unique under
the condition that the principal symbol be preserved at each order.
\end{thm}

The singular values (\ref{Reson}) of the shift $\d=\m-\l$ are called resonances
and lead to special and interesting modules.

Theorem \ref{confEquivIsomTh} is particular case of a more general result
recently obtained in \cite{DLO}.

\begin{thm}\label{confEquivIsomThRes}
For each resonant value of $\delta$, there exist particular pairs
$(\l,\m)$ of weights such that the $\so(p+1,q+1)$-modules $\cS^2_\delta$ and
$\cD^2_{\l,\m}$ are isomorphic, namely
\begin{equation}
\setlength{\extrarowheight}{8pt}
\begin{array}{|c||c|c|c|c|c|}
\hline
\delta & \frac{2}{n} & \frac{n+2}{2n} & 1 & \frac{n+1}{n} &
\frac{n+2}{n}\\[8pt]
\hline
\hline
\lambda &
\frac{n-2}{2n} &
0,\frac{n-2}{2n} & 0 & 0, -\frac{1}{n} & -\frac{1}{n}\\[8pt]
%&&&&&\\
\hline
\mu & \frac{n+2}{2n} & \frac{n+2}{2n},1 & 1 &
\frac{n+1}{n},1 & \frac{n+1}{n}\\[8pt]
%&&&&&\\
\hline
\end{array}
\label{TheArray}
\end{equation}
\end{thm}

We will show 
in Section \ref{ProofResonThm} 
that the isomorphism is not unique;
there exists, actually, a one-parameter family of such isomorphisms in each resonant
case.

\begin{rmk}{\rm
This point of view about equivariant quantization was adopted in~\cite{LO} where a
projectively-equivariant symbol calculus and quantization was introduced if $M$ is
endowed with a flat projective structure. In this case the group of (local)
symmetries is $G=\SL(n+1,\bbR)$ with
$n=\dim(M)$. See also \cite{Lec} for a cohomological treatment of this subject.
Bearing in mind that the best-known geometries associated with a local and
maximal (see Section~\ref{ConclusionOutlook}) symmetry group are the projective and
conformal geometries, we have been led to look, in the same spirit, for a
conformally equivariant quantization.
}
\end{rmk}

\begin{rmk}{\rm
In the particular case $n=1$, the projective and conformal symmetries coincide; our
results are in full accordance with those obtained in \cite{GO,Gar,CMZ} and the
resonances are simply $\{1,\frac{3}{2},2\}$. 
}
\end{rmk}

\goodbreak

%%%%%%%%%%%%%%%%%%%%%%%%%%%%%%%%%%%%%%%%%%%%%%%%%%%%%%%%%%%%%%%%%%%%%%%%%%%%%%%%%%%%
%\subsection{Conformal invariance}
%%%%%%%%%%%%%%%%%%%%%%%%%%%%%%%%%%%%%%%%%%%%%%%%%%%%%%%%%%%%%%%%%%%%%%%%%%%%%%%%%%%%

It turns out that this isomorphism $\cQ_{\l,\m}$ makes sense for an arbitrary
pseudo-Riemannian manifold (not necessarily conformally flat). The fundamental
property of this isomorphism is that it depends only on the conformal class $[\rg]$
of the metric (i.e., it is conformally invariant) --- see Theorem \ref{invConfQ2}.

We have no precise uniqueness theorem here, but we will show that the condition
of conformal invariance uniquely determines $\cQ_{\l,\m}$ in some natural class
of differential linear maps. This enables us to introduce a conformally invariant
quantization on the cotangent bundle of a pseudo-Riemannian manifold.
Note that we understand the term ``quantization'' in a somewhat generalized sense
as~$\l$ and~$\m$ remain essentially arbitrary.
In the case $\l=\m=\half$, we recover the usual
terminology using the Hilbert space of half-densities considered in the framework of
geometric quantization.

As an illustration of our general results, we consider a number of examples. The
celebrated Yamabe-Laplace operator is, among others, naturally  included into our
considerations. It appears as the quantized geodesic flow in one of the resonant
cases~(\ref{TheArray}); the corresponding symbol is itself  conformally invariant.
This explains why the Yamabe-Laplace operator is the unique conformally invariant
Laplace-Beltrami operator. It should be stressed, however, that it is unjustified to
consider the Yamabe-Laplace operator as quantum Hamiltonian for the geodesic flow.

%%%%%%%%%%%%%%%%%%%%%%%%%%%%%%%%%%%%%%%%%%%%%%%%%%%%%%%%%%%%%%%%%%%%%%%%%%%%%%%%%%%%
%\subsection{Organization of the paper}
%%%%%%%%%%%%%%%%%%%%%%%%%%%%%%%%%%%%%%%%%%%%%%%%%%%%%%%%%%%%%%%%%%%%%%%%%%%%%%%%%%%%

The paper is organized as follows.

In Section \ref{BasicDefinitions} we recall the basic definitions
concerning the space of differential operators on tensor densities, as well as the
space of symbols. We put emphasis on their $\Diff(M)$- and $\Vect(M)$-module
structures.

We present, in Section \ref{explicitFormulae}, an explicit intrinsic
formula for the isomorphism~$\cQ_{\l,\m}$ which defines our conformal\-ly
equivariant quantization map $\cQ_{\l,\m;\hbar}$ for an arbitrary
pseudo-Riemannian manifold. We also prove the conformal invariance
of~$\cQ_{\l,\m}$. 

Section \ref{ApplicationSection} provides specific
examples, namely the quantization of the geodesic flow, and of the magnetic minimal
coupling prescription. 
We also give examples of quantized Hamiltonians pertaining
to the resonant cases. 

We develop in Section \ref{WeylBrauerSection} the algebraic theory of
Euclidean invariants which we use in the proofs of the uniqueness
theorem for the conformally equivariant quantization map.

In the last Sections \ref{ConfEquivSection} and \ref{proofsMainResults} we give
the technical proofs of the main theorems.

%%%%%%%%%%%%%%%%%%%%%%%%%%%%%%%%%%%%%%%%%%%%%%%%%%%%%%%%%%%%%%%%%%%%%%%%%%%%%%%%%%%%
%%%%%%%%%%%%%%%%%%%%%%%%%%%%%%%%%%%%%%%%%%%%%%%%%%%%%%%%%%%%%%%%%%%%%%%%%%%%%%%%%%%%
\section{Differential operators and symbols}
\label{BasicDefinitions}
%%%%%%%%%%%%%%%%%%%%%%%%%%%%%%%%%%%%%%%%%%%%%%%%%%%%%%%%%%%%%%%%%%%%%%%%%%%%%%%%%%%%
%%%%%%%%%%%%%%%%%%%%%%%%%%%%%%%%%%%%%%%%%%%%%%%%%%%%%%%%%%%%%%%%%%%%%%%%%%%%%%%%%%%%

%%%%%%%%%%%%%%%%%%%%%%%%%%%%%%%%%%%%%%%%%%%%%%%%%%%%%%%%%%%%%%%%%%%%%%%%%%%%%%%%%%%%
\subsection{Differential operators on tensor densities}
%%%%%%%%%%%%%%%%%%%%%%%%%%%%%%%%%%%%%%%%%%%%%%%%%%%%%%%%%%%%%%%%%%%%%%%%%%%%%%%%%%%%

Let us start with the definition of the $\Diff(M)$-module $\cD_{\l,\m}(M)$ (or 
$\cD_{\l,\m}$ for short) of differential operators on a smooth manifold $M$ with
$\l,\m\in\bbR\,(\hbox{or\ }\bbC)$.

Consider the determinant bundle $\Lambda^nT^*M\rightarrow M$.
Let us recall that a tensor density of degree $\l$ on $M$ is a
smooth section, $\phi$, of the line bundle
$\mod{\Lambda^nT^*M}^{\otimes\l}$. The space of tensor densities of
degree $\l$ is naturally a $\Diff(M)$-module which we call~$\cF_\l$.

It is evident that $\cF_0=C^{\infty}(M)$; if $M$ is oriented, the space $\cF_1$
coincides with the space of differential $n$-forms: $\cF_1=\Omega^n(M)$. 

\begin{defi} 
An operator $A:\cF_\l\to\cF_\m$ is called a local operator on~$M$ if
for all $\phi\in\cF_\l$ one has $\Supp(A(\phi))\subset\Supp(\phi)$.
\end{defi}

It is a classical result (see \cite{Pee}) that such operators are in fact
locally given by differential operators. The space $\cD_{\l,\m}$ of differential
operators from $\l$-densities to $\m$-densities on $M$ is  naturally a
$\Diff(M)$-module.

There is a filtration
$\cD^0_{\l,\m}\subset\cD^1_{\l,\m}\subset\cdots\subset\cD^k_{\l,\m}\subset\cdots$,
where the module of zero-order operators $\cD^0_{\l,\m}\cong\cF_{\m-\l}$
consists of multiplication by $(\m-\l)$-densities. The higher-order modules are
defined by induction: $A\in\cD^k_{\l,\m}$ if $[A,f]\in\cD^{k-1}_{\l,\m}$ for
every $f\in{}C^\infty(M)$.

To our knowledge, the whole family of modules of differential operators viewed as
a deformation were first studied in~\cite{DO} in the case $\l=\m$. (See also
\cite{LMT,LO,GO,Gar,Mat}.) 

%%%%%%%%%%%%%%%%%%%%%%%%%%%%%%%%%%%%%%%%%%%%%%%%%%%%%%%%%%%%%%%%%%%%%%%%%%%%%%%%%%%%
\subsection{Classical examples}
%%%%%%%%%%%%%%%%%%%%%%%%%%%%%%%%%%%%%%%%%%%%%%%%%%%%%%%%%%%%%%%%%%%%%%%%%%%%%%%%%%%%

\indent\indent 
(a) The best known example is the Sturm-Liouville operator
$L=(d/dx)^2+u(x)$ in the one-dimensional case, $M=S^1$. It should, indeed, be
considered as an element
$L\in\cD^2_{-\frac{1}{2},\frac{3}{2}}$ as $\l=-1/2$ and
$\m=3/2$ are the only degrees for which its form is preserved by the
action of $\Diff(S^1)$.

(b) Again, in the one-dimensional case, the study of
the modules $\cD^k_{\frac{1-k}{2},\frac{1+k}{2}}$ goes back to the pioneering
work of Wilczynski~\cite{Wil}.

(c) Yet another remarkable example is provided by the Yamabe-Laplace operator
$A=\Delta-(n-2)/(4(n-1))\,R$, where $\Delta$ is the usual
Laplace-Beltrami operator and~$R$ the scalar curvature on a
(pseudo-)Riemannian mani\-fold $(M,\rg)$ of dimension $n\geq2$. (See,
e.g.~\cite{Bes}.) This operator has been extensively used in the mathematical
and physical literature because of its characteristic property of being
invariant under conformal changes of metrics. It is well known that
$A\in\cD^2_{\frac{n-2}{2n},\frac{n+2}{2n}}$. 

Observe that, for $n=1$, the latter module of differential operators precisely
coincides with the Sturm-Liouville module. We will see that
this is by no means accidental and will prove below (Section
\ref{SturmLiouvilleSect}) that the suitably regularized Yamabe operator equals
$\Delta-S(\varphi)/(2\rg)$, where $S$ is the Schwarzian derivative and $\varphi$
the diffeomorphism which defines the metric $\rg=\varphi^*(dx^2)$.

(d) The special module $\cD_{\half,\half}$ has been introduced in the context of
geometric quantization by Blattner \cite{Bla} and Kostant \cite{Kos}. This module
will also naturally arise in our quantization procedure.

%%%%%%%%%%%%%%%%%%%%%%%%%%%%%%%%%%%%%%%%%%%%%%%%%%%%%%%%%%%%%%%%%%%%%%%%%%%%%%%%%%%%
\subsection{The modules $\cF_\l$ and $\cD_{\l,\m}$}
%%%%%%%%%%%%%%%%%%%%%%%%%%%%%%%%%%%%%%%%%%%%%%%%%%%%%%%%%%%%%%%%%%%%%%%%%%%%%%%%%%%%

If $M$ is orientable, which we will assume throughout the paper, then $\cF_\l$ can be
identified with
$C^{\infty}(M)$ as a vector space. Given a volume form, $\vol$, on~$M$, one can
write any $\lambda$-tensor density as $\phi=f\mod{\vol}^\l$ with
$f\in{}C^\infty(M)$, and define the $\Diff(M)$-module structure of
$\cF_{\lambda}$ via the action of $\varphi\in\Diff(M)$:
\begin{equation}
\varphi_\lambda(f) =
\varphi_*(f)\,\left|\frac{\varphi_*\vol}{\vol}\right|^\lambda.
\label{Faction}
\end{equation}

With this identification, the module $\cD_{\l,\m}$ can be viewed as a
two-parameter family of the standard module $\cD_{0,0}$ of differential
operators on smooth functions~$\cF_0$. The natural $\Diff(M)$-action on
$\cD_{\l,\m}$ then reads
\begin{equation}
\varphi_{\l,\m}(A) = 
\varphi_\m\circ{}A\circ{}\varphi^{-1}_\l.
\label{Daction}
\end{equation}

The expression of a differential operator $A\in\cD^k_{\l,\m}$ in a local
coordinate system~$(x^i)$ is then
\begin{equation}
A
=
A_k^{{i_1}\ldots{i_k}}\partial_{i_1}\ldots\partial_{i_k} 
+
\cdots 
+
A_1^i\partial_i 
+
A_0
\label{DiffOp}
\end{equation}
where $\partial_i=\partial/\partial x^i$ and 
$A_\ell^{{i_1}\ldots{i_\ell}}\in{}C^{\infty}(M)$
with $\ell=0,1,\ldots,k$. 
(From now on we suppose a summation over repeated indices.)

The infinitesimal version of the action (\ref{Daction}) is
\begin{equation}
L^{\l,\m}_X(A) = 
L^{\m}_X\,A - A\,L^{\l}_X
\label{Dinf}
\end{equation}
where $X\in\Vect(M)$, while the infinitesimal version of the action
(\ref{Faction}) is given by the Lie derivative on $\cF_\l$, namely
\begin{equation}
L^\l_X(f) = 
X(f)+\l\,\Div(X)\,f.
\label{Finf}
\end{equation}

\goodbreak

%%%%%%%%%%%%%%%%%%%%%%%%%%%%%%%%%%%%%%%%%%%%%%%%%%%%%%%%%%%%%%%%%%%%%%%%%%%%%%%%%%%%
\subsection{The module of symbols $\cS_\d$}
%%%%%%%%%%%%%%%%%%%%%%%%%%%%%%%%%%%%%%%%%%%%%%%%%%%%%%%%%%%%%%%%%%%%%%%%%%%%%%%%%%%%

Consider the space $\cS=\Gamma(S(TM))$ of contravariant symmetric tensor
fields on~$M$ which is naturally a $\Diff(M)$-module. We can locally identify $\cS$
with the space of polynomials
\begin{equation}
P(\xi)=\sum_{\ell=0}^k{P_\ell^{i_1\ldots{}i_\ell}\xi_{i_1}\cdots\xi_{i_\ell}},
\label{Symb}
\end{equation}
with
$P_\ell^{i_1\ldots i_\ell}\in C^\infty(M)$, on the cotangent bundle of
$M$.

\begin{defi}
\label{defSymb}
The one-parameter family of $\Diff(M)$-actions on $\cS$:
\begin{equation}
\varphi_\d(P)=
\varphi_*(P)\left|\frac{\varphi_*\vol}{\vol}\right|^\d
\label{DiffActionSymb}
\end{equation}
identifies the space $\cS$ with the $\Diff(M)$-module
$\cS\otimes\cF_\d$. We denote this module by~$\cS_\d$.
\end{defi}

We will need in the sequel the infinitesimal version of the $\Diff(M)$-action on
$\cS_\d$. 
The action of $\Vect(M)$ on $\cS_\d$ deduced from (\ref{DiffActionSymb}) reads:
\begin{equation}
L_X^\d(P)=
L_X(P)+\d\,\Div(X)\,P
\label{VectActionSymb}
\end{equation}
where
\begin{equation}
L_X=
X^i\frac{\partial}{\partial x^i}
-
\xi_j\partial_iX^j\frac{\partial}{\partial \xi_i}
\label{LieDer}
\end{equation}
is the cotangent lift of $X\in\Vect(M)$.

\goodbreak

Again, there is a filtration
$\cS^0_\d\subset\cS^1_\d\subset\cdots\subset\cS^k_\d\subset\cdots$,
where $\cS^k_\d$ denotes the space of symbols of degree less or equal to $k$. 
In contrast to the filtration on the space $\cD_{\l,\m}$ of differential operators, the
above filtration on the space of symbols actually leads to a
$\Diff(M)$-invariant graduation 
\begin{equation}
\cS_\d=\bigoplus_{k=0}^\infty{\cS_{k,\d}}
\label{grad}
\end{equation}
where $\cS_{k,\d}$ denotes the space of homogeneous polynomials
(isomorphic to $\cS^k_\d/\cS^{k-1}_\d$).

%%%%%%%%%%%%%%%%%%%%%%%%%%%%%%%%%%%%%%%%%%%%%%%%%%%%%%%%%%%%%%%%%%%%%%%%%%%%%%%%%%%%
%%%%%%%%%%%%%%%%%%%%%%%%%%%%%%%%%%%%%%%%%%%%%%%%%%%%%%%%%%%%%%%%%%%%%%%%%%%%%%%%%%%%
\section{Explicit formul{\ae} for the quantization map}
\label{explicitFormulae}
%%%%%%%%%%%%%%%%%%%%%%%%%%%%%%%%%%%%%%%%%%%%%%%%%%%%%%%%%%%%%%%%%%%%%%%%%%%%%%%%%%%%
%%%%%%%%%%%%%%%%%%%%%%%%%%%%%%%%%%%%%%%%%%%%%%%%%%%%%%%%%%%%%%%%%%%%%%%%%%%%%%%%%%%%

There is no fully $\Diff(M)$-equivariant quantization since the modules
$\cD_{\l,\m}$ are not isomorphic to the module $\cS_{\m-\l}$ of symbols. One is
thus led to impose some extra geometric structure on $M$ and to look for a
symbol calculus, equivariant with respect to the automorphisms of this
structure. 

In this article, we will assume---unless otherwise stated---that the manifold $M$
is endowed with a flat conformal structure.

%%%%%%%%%%%%%%%%%%%%%%%%%%%%%%%%%%%%%%%%%%%%%%%%%%%%%%%%%%%%%%%%%%%%%%%%%%%%%%%%%%%%
\subsection{A compendium on conformally flat structures}\label{ConfFlatStruct}
%%%%%%%%%%%%%%%%%%%%%%%%%%%%%%%%%%%%%%%%%%%%%%%%%%%%%%%%%%%%%%%%%%%%%%%%%%%%%%%%%%%%

A conformal structure on a manifold $M$ is given by a smooth field $[\,g\,]$
of directions of metrics. This structure is called flat if $M$ can be
locally identified with~$\bbR^n$ endowed with the canonical action of the conformal
Lie algebra $\so(p+1,q+1)$, where $n=p+q$. 

The Lie algebra $\so(p+1,q+1)\subset\Vect(\bbR^n)$ is generated by the
vector fields:
\begin{equation}
%\left\{
\matrix{
X_i   &=&\displaystyle \frac{\partial}{\partial x^i}\;,\hfill\cr
\noalign{\smallskip}
X_{ij}&=&\displaystyle x_i\frac{\partial}{\partial x^j}-
x_j\frac{\partial}{\partial x^i}\;,\hfill\cr
\noalign{\smallskip}
X_0   &=&\displaystyle x^i\frac{\partial}{\partial x^i}\;,\hfill\cr
\noalign{\smallskip}
\bar X_i &=&\displaystyle x_jx^j\frac{\partial}{\partial x^i}-
2x_ix^j\frac{\partial}{\partial x^j}\hfill\cr
}
%\right.
\label{ConfAlg}
\end{equation}
with $i,j=1,\ldots,n$; we have used the notation $x_i=g_{ij}x^j$ where the flat
metric $g=\diag(1,\ldots,1,-1,\ldots,-1)$ has trace $p-q$.

\goodbreak

The subalgebra generated by the vector fields $X_i$ and $X_{ij}$ is the
Euclidean Lie algebra $\se(p,q)=\so(p,q)\ltimes\bbR^n$. The operator $X_0$ is the
generator of homotheties while the vector fields $\bar X_i$ generate
inversions.

\begin{rmk}{\rm
(a) It is well known that the conformal flatness of a
$n$-dimensional pseudo-Riemannian mani\-fold is equivalent to the vanishing of the
Weyl curvature tensor if $n\geq4$, and to that of the Weyl-Schouten curvature tensor
if $n=3$ \cite{Bes}. All two-dimensional pseudo-Riemannian mani\-folds are
conformally flat.

(b) In the one-dimensional case the conformal Lie algebra is
isomorphic to the projective Lie algebra since $\so(2,1)\cong\Sl(2,\bbR)$. 
}
\end{rmk}

\goodbreak

%%%%%%%%%%%%%%%%%%%%%%%%%%%%%%%%%%%%%%%%%%%%%%%%%%%%%%%%%%%%%%%%%%%%%%%%%%%%%%%%%%%%
\subsection{Conformal equivariance}
%%%%%%%%%%%%%%%%%%%%%%%%%%%%%%%%%%%%%%%%%%%%%%%%%%%%%%%%%%%%%%%%%%%%%%%%%%%%%%%%%%%%

Let $M$ be endowed with a flat conformal structure: there exists a local
action of the group $\Or(p+1,q+1)$ on $M$, which enables us to restrict the
$\Diff(M)$-modules~$\cD_{\l,\m}$ to the conformal group. Our problem
amounts then to the determination of intertwining differentiable linear maps
$\cQ_{\l,\m}^k$
between the $\so(p+1,q+1)$-modules $\cS^k_{\l-\m}$ and $\cD^k_{\m,\l}$.

\goodbreak

Here, we give the solution for the case $k=2$ which is the most relevant
one for applications. Indeed, the existence and uniqueness of a conformally
equivariant quantization map for any order $k$ has recently been  established in
\cite{DLO}; however, no explicit formula is available.

%%%%%%%%%%%%%%%%%%%%%%%%%%%%%%%%%%%%%%%%%%%%%%%%%%%%%%%%%%%%%%%%%%%%%%%%%%%%%%%%%%%%
\subsection{Expression in adapted coordinates}
%%%%%%%%%%%%%%%%%%%%%%%%%%%%%%%%%%%%%%%%%%%%%%%%%%%%%%%%%%%%%%%%%%%%%%%%%%%%%%%%%%%%

If we fix the local coordinate system on $M$ for which the generators of
$\so(p+1,q+1)$ retain the form (\ref{ConfAlg}), we have the following  

\begin{thm}\label{ThmQ2}
For any dimension, $n$, and any $\d$ as in (\ref{Reson}) the unique conformally
equivariant isomorphism
$\cQ_{\l,\m}:\cS_\d^2\to\cD_{\l,\m}^2$, viz 
$$
P\longmapsto\cQ_{\l,\m}(P)=
A_2^{ij}\partial_i\partial_j+A_1^i\partial_i+A_0
$$
that preserves the principal symbol is as
follows
\begin{equation}
\left\{\begin{array}{lcl}
A_2^{ij}  
&=& 
P_2^{ij}\hfill\\[4pt]
A_1^i  
&=& 
P_1^i+\beta_1\,\partial_jP_2^{ij}+
\beta_2\,g^{ij}g_{k\ell}\,\partial_j\,P_2^{k\ell}\\[4pt]
A_0  
&=& 
P_0+\alpha\,\partial_iP_1^i+
\beta_3\,\partial_{ij}P_2^{ij}+
\beta_4\,g^{ij}g_{k\ell}\,\partial_{ij}P_2^{k\ell}
\end{array}
\right.
\label{Symbol}
\end{equation}
where $P(\xi)
=
P_2^{ij}\xi_i\xi_j+P_1^i\xi_i+P_0
\in\cS^2_\d$; the numerical coefficients are given by
\begin{equation}
\a=\frac{\l}{1-\d}
\label{alpha}
\end{equation}
and
\begin{equation}
\begin{array}{l}
\beta_1
=
\displaystyle 
\frac{2(n\l+1)}{2+n(1-\d)}\hfill\\[12pt]
\beta_2
=
\displaystyle 
\frac{n(\l+\m-1)}{(2+n(1-\d))(2-n\d)}\hfill\\[12pt]
\beta_3
=
\displaystyle 
\frac{n\l(n\l+1)}{(1+n(1-\d))(2+n(1-\d))}\hfill\\[12pt]
\beta_4
=
\displaystyle 
\frac{n\l(n^2\m(2-\l-\m)+2(n\l+1)^2-n(n+1))}%
{(1+n(1-\d))(2+n(1-\d))(2+n(1-2\d))(2-n\d)}.
\end{array}
\label{beta14}
\end{equation}
\end{thm}
We will prove this theorem in Section \ref{proofsMainResults}.

In the one-dimensional case, $n=1$, this formula can
be written as
\begin{equation}
\left\{
\begin{array}{lcl}
A_2  
&=& 
P_2\hfill\\[10pt]
A_1  
&=& 
P_1 
+ 
\displaystyle\frac{2\l+1}{2-\d}\,P'_2
\hfill\\[10pt]
A_0  
&=& 
P_0
+
\displaystyle\frac{\l}{1-\d}\,P'_1
+
\displaystyle\frac{\l(2\l+1)}{(3-2\d)(2-\d)}\,P''_2.
\hfill%\\[10pt]
\end{array}
\right.
\label{oneDimQ}
\end{equation}
\begin{rmk}{\rm
We record that the formula (\ref{Symbol}) will be the main ingredient of the
proof of Proposition
\ref{ProQ1} and Theorem \ref{ThmQ2Int} below establishing intrinsic formul{\ae} for
the quantization map.
}
\end{rmk}

\begin{rmk}{\rm
The projectively equivariant symbol map (and its inverse, the
quantization map) has been constructed in~\cite{LO} in the special case
$\l=\m$ in any dimension. (See also \cite{Lec} for arbitrary $\l$ and~$\m$, and
\cite{CMZ,Gar} for the one-dimensional case.) 
}
\end{rmk}

%%%%%%%%%%%%%%%%%%%%%%%%%%%%%%%%%%%%%%%%%%%%%%%%%%%%%%%%%%%%%%%%%%%%%%%%%%%%%%%%%%%%
\subsection{The covariant derivative of densities}
%%%%%%%%%%%%%%%%%%%%%%%%%%%%%%%%%%%%%%%%%%%%%%%%%%%%%%%%%%%%%%%%%%%%%%%%%%%%%%%%%%%%

Given a conformally flat manifold $M$ of signature $p-q$, one can  choose, locally, a
pseudo-Riemannian metric $\rg$ which represents the conformal class of the manifold.
We will denote by $\nabla$ the Levi-Civita connection.
Let us now recall the definition of the covariant derivative of densities. If
$\phi\in\cF_\l$, then $\nabla\phi\in\Omega^1(M)\otimes\cF_\l$ is defined by
$\nabla\phi=df\otimes\mod{\vol}^\l$, using the local representation
$\phi=f\mod{\vol}^\l$ with $f\in{}C^\infty(M)$.

Choose an arbitrary coordinate system $(x^i)$ on $M$ (with associated
coordinate system $(\xi_i,x^i)$ on $T^*M$); one has, for every $\phi\in\cF_\l$, the
local expression
\begin{equation} 
\nabla_i\phi=\partial_i\phi-\l\Gamma_i\phi
\label{nablaphi}
\end{equation} 
with $\Gamma_i=\Gamma_{ij}^j$.

\goodbreak

%%%%%%%%%%%%%%%%%%%%%%%%%%%%%%%%%%%%%%%%%%%%%%%%%%%%%%%%%%%%%%%%%%%%%%%%%%%%%%%%%%%%
\subsection{Intrinsic expression for first-order polynomials}
%%%%%%%%%%%%%%%%%%%%%%%%%%%%%%%%%%%%%%%%%%%%%%%%%%%%%%%%%%%%%%%%%%%%%%%%%%%%%%%%%%%%

The intrinsic form of the isomorphism (\ref{confEquivIsom}) for
first-order polynomials is given by
\begin{pro}\label{ProQ1}
For any $\d\neq1$ the unique conformally equivariant isomorphism
$\cQ_{\l,\m}:\cS_\d^1\to\cD_{\l,\m}^1$ that preserves the principal symbol is given
by
\begin{equation}
\cQ_{\l,\m}(P)=
P_1^i\nabla_i+ \a\nabla_i(P_1^i) + P_0
\label{Q1}
\end{equation}
where $P(\xi)=P_1^i\xi_i+P_0$ and $\alpha$ as in (\ref{alpha}).
\end{pro}

\goodbreak

It can be verified that $\cQ_{\l,\m}$ in (\ref{Q1}) is, actually, equivariant with
respect to the full Lie algebra $\Vect(M)$. This formula holds in any dimension. 

\goodbreak

\begin{rmk}{\rm
In the resonant case, $\d=1$, the modules are still isomorphic if and
only if $(\l,\m)=(0,1)$ as given by Theorem \ref{confEquivIsomThRes}. The isomorphism
is not unique and given by the formula (\ref{Q1}) with arbitrary $\a$.
}
\end{rmk}

%%%%%%%%%%%%%%%%%%%%%%%%%%%%%%%%%%%%%%%%%%%%%%%%%%%%%%%%%%%%%%%%%%%%%%%%%%%%%%%%%%%%
\subsection{Intrinsic expression for quadratic polynomials in the multi-dimensional
case}
%%%%%%%%%%%%%%%%%%%%%%%%%%%%%%%%%%%%%%%%%%%%%%%%%%%%%%%%%%%%%%%%%%%%%%%%%%%%%%%%%%%%

Let us now give the explicit expression of the isomorphism (\ref{confEquivIsom}) for
homogeneous second-order symbols in the case $\dim(M)\geq3$.

\begin{thm}\label{ThmQ2Int}
If $n\geq3$, for any $\d$ as in (\ref{Reson}) the unique conformally
equivariant isomorphism (\ref{confEquivIsom}) that preserves the principal symbol is
as follows
\begin{equation}
\begin{array}{lcl}
\cQ_{\l,\m}(P)  
&=& 
P^{ij}\nabla_i\nabla_j \\[10pt]
&&+ \left(
\beta_1\nabla_iP^{ij}+\beta_2\,\rg^{ij}\rg_{k\ell}\nabla_iP^{k\ell}
\right)\nabla_j
\\[10pt]
&& 
+\beta_3\nabla_i\nabla_j(P^{ij})
+\beta_4\,\rg^{ij}\rg_{k\ell}\nabla_i\nabla_j(P^{k\ell})
+\beta_5R_{ij}P^{ij}
+\beta_6R\,\rg_{ij}P^{ij}\\%[10pt]
\end{array}
\label{Q2}
\end{equation}
where $P(\xi)=P^{ij}\xi_i\xi_j$; the coefficients $\beta_1,\ldots,\beta_4$ are
given by (\ref{beta14}) and
\begin{equation}
\begin{array}{l}
\beta_5
=
\displaystyle 
\frac{n^2\l(\m-1)}{(n-2)(1+n(1-\d))}\\[10pt]
\beta_6
=
\displaystyle 
\frac{n^2\l(\m-1)(n\d-2)}{(n-1)(n-2)(1+n(1-\d))(2+n(1-2\d))}
\end{array}
\label{beta56}
\end{equation}
and $R_{ij}$ (resp. $R$) denote the Ricci tensor components (resp. the scalar
curvature) of the metric $\rg$.
\end{thm}
This theorem will be proven in Section \ref{proofsMainResults}; the main ingredient
of its proof will be Theorem \ref{ThmQ2}. 

We also use, in the sequel, the notation $\Ric=R_{ij}\,dx^i\otimes{}dx^j$
for the Ricci tensor.

\begin{rmk}{\rm Another quantization formula for second-order polynomials has been
proposed in \cite{LQ} using a (pseudo-)Riemannian metric on $M$ and the local
identification of $T^*M$ with $\bbR^{2n}$ endowed with its standard $\ssp(2n,\bbR)$
action. }
\end{rmk}

%%%%%%%%%%%%%%%%%%%%%%%%%%%%%%%%%%%%%%%%%%%%%%%%%%%%%%%%%%%%%%%%%%%%%%%%%%%%%%%%%%%%
\subsection{Lower-dimensional cases and Schwarzian derivatives}
%%%%%%%%%%%%%%%%%%%%%%%%%%%%%%%%%%%%%%%%%%%%%%%%%%%%%%%%%%%%%%%%%%%%%%%%%%%%%%%%%%%%

The general formula (\ref{Q2}) for the conformally equivariant map is
obviously non applicable in the cases $n=1$ and $n=2$. We must therefore consider
each of these cases separately.

Let us start with the one-dimensional case for which all metrics are equivalent.
In this case, $M=S^1$ say, the metric retains the form $\rg=\varphi^*(dx^2)$ for some
$\varphi\in\Diff(S^1)$ with $x$ an arbitrary coordinate. 
\begin{thm}
\label{ThmQ2n1}
If $n=1$, and $\d\neq\frac{3}{2},2$,
the unique conformally
equivariant isomorphism (\ref{confEquivIsom}) that preserves the principal symbol is
given by
\begin{equation}
\cQ_{\l,\m}(P)
=
P\,\nabla^2
+
\frac{2\l+1}{2-\d}\,(\nabla{P})\nabla
+
\frac{\l(2\l+1)}{(3-2\d)(2-\d)}\,(\nabla^2{P})
-
\,\frac{2\l(\m-1)}{3-2\d}\,\frac{S(\varphi)}{\rg}
\label{Q2One}
\end{equation}
where $P(\xi)=P\xi^2$ and
\begin{equation}
S(\varphi)
=
\frac{\varphi'''}{\varphi'}
-
\frac{3}{2}\left(\frac{\varphi''}{\varphi'}\right)^2
\label{Schwarzian}
\end{equation}
is the Schwarzian derivative of $\varphi$.
%The Laplace operator takes the form 
%\begin{equation}
%\Delta_\rg
%=
%\frac{1}{\rg}\left(
%\frac{d^2}{dx^2} - \frac{\varphi''}{\varphi'}\frac{d}{dx}
%\right).
%\label{LaplaceOne}
%\end{equation}
\end{thm}
Comparison with the expression (\ref{Q2}) strengthens
the saying according to which the Schwarzian derivative is nothing but
``curvature''.

The two-dimensional case, $n=2$, is especially interesting since all
surfaces~$(M,\rg)$ are conformally flat. 
The Riemann uniformization theorem can be invoked to express the metric (locally)
as
\begin{equation}
\rg=F^{-1}\varphi^*\rg_0
\label{uniformization}
\end{equation}
where $\varphi$ is a conformal diffeomorphism of $M$, and
$F\in{}C^\infty(M,\bbR^*_+)$, and
$\rg_0$ is a metric of constant curvature. Let us emphasize that this weaker
form of the uniformization theorem still holds in the Lorentz case
(see, e.g., \cite{Wei}).

There exists in the recent literature an interesting generalization of
the Schwar\-zian derivative for conformal diffeomorphisms in the multi-dimensional
case. In the situation (\ref{uniformization}) with $F=e^{2f}$, the Schwarzian
derivative of $\varphi$ is defined \cite{OS,Car} as the  symmetric
twice-covariant tensor $S(\varphi)$ such that
\begin{equation}
S(\varphi)(X,Y)=
X(Yf)-(\nabla_XY)f-(Xf)(Yf)+\half\Vert{df}\Vert^2_\rg\,\rg(X,Y)
\label{Carne}
\end{equation}
for any $X,Y\in\Vect(M)$. In our notation, it reads
\begin{equation}
S(\varphi)=
\frac{1}{2F}\nabla{dF}
-
\frac{3}{4F^2}\,dF\otimes{}dF
+
\frac{1}{8F^2}\,\rg^{-1}(dF,dF)\,\rg.
\label{multiDimSchwarzian}
\end{equation}
This new object will enter naturally the expression of the conformally
equivariant map (\ref{confEquivIsom}) for surfaces.

\goodbreak

Note that the definition (\ref{multiDimSchwarzian}) yields the classical Schwarzian
derivative in the one-dimensional case. 

\goodbreak

\begin{thm}\label{ThmQ2n2}
If $n=2$, for any $\d$ as in (\ref{Reson}) the unique conformally
equivariant isomorphism (\ref{confEquivIsom}) preserving the principal symbol reads
\begin{equation}
\begin{array}{lcl}
\cQ_{\l,\m}(P)  
&=& 
P^{ij}\nabla_i\nabla_j \\[10pt]
&&+ \left(
\beta_1\nabla_iP^{ij}+\beta_2\,\rg^{ij}\rg_{k\ell}\nabla_iP^{k\ell}
\right)\nabla_j
\\[10pt]
&& 
+\beta_3\nabla_i\nabla_j(P^{ij})
+\beta_4\,\rg^{ij}\rg_{k\ell}\nabla_i\nabla_j(P^{k\ell}) \\[10pt]
&&
+\displaystyle\frac{4\l(\m-1)}{2\d-3}\Big(S(\varphi)_{ij}P^{ij}
+\frac{1}{8(\d-1)}R\,\rg_{ij}P^{ij}\Big)\\%[10pt]
\end{array}
\label{Q2S}
\end{equation}
where $P(\xi)=P^{ij}\xi_i\xi_j$ and $S(\varphi)$ is as in
(\ref{multiDimSchwarzian}) while $R$ denotes the scalar curvature of $\rg$; the
coefficients $\beta_1,\ldots,\beta_4$ are given by (\ref{beta14}).
\end{thm}
Notice that the scalar curvature is related
to the trace of the  Schwarzian derivative by
\begin{equation}
R = -2\rg^{ij}S(\varphi)_{ij}
\label{TraceSchwarz}
\end{equation}
provided (\ref{uniformization}) holds.

We defer the proofs of Theorems \ref{ThmQ2n1} and \ref{ThmQ2n2}
to Section
\ref{proofsMainResults}.

%%%%%%%%%%%%%%%%%%%%%%%%%%%%%%%%%%%%%%%%%%%%%%%%%%%%%%%%%%%%%%%%%%%%%%%%%%%%%%%%%%%%
\subsection{Conformal invariance}
%%%%%%%%%%%%%%%%%%%%%%%%%%%%%%%%%%%%%%%%%%%%%%%%%%%%%%%%%%%%%%%%%%%%%%%%%%%%%%%%%%%%

The preceding formul{\ae} 
are intrinsic and therefore extend to any pseudo-Rieman\-nian manifold $(M,\rg)$,
not necessarily conformally flat. They altogether define the
map~(\ref{confEquivIsom}).
\begin{thm}
The map $\cQ_{\l,\m}:\cS^2_{\m-\l}\longrightarrow\cD^2_{\l,\m}$ defined by
(\ref{Q1}), (\ref{Q2}), (\ref{Q2One}) and (\ref{Q2S}) is conformally
invariant,i.e., it depends only on the conformal class of the metric.
\label{invConfQ2}
\end{thm}
\begin{proof}
Let us choose another metric $\widehat\rg=F\rg$ with $F$ a strictly positive valued
function.  In the special case of conformally flat manifolds, the map
$\cQ_{\l,\m}$ is given, in an adapted coordinate system, by
Theorem~\ref{ThmQ2}. Now, the adapted coordinate systems for~$\rg$
and~$\widehat\rg$ are the same. This proves the theorem in the conformally flat
case, in particular for~$n=1$ and $n=2$ in full generality. 

\goodbreak

The case of an arbitrary pseudo-Rieman\-nian manifold needs a separate proof
which goes as follows. We have
\begin{equation}
\widehat\Gamma_{ij}^k
=
\Gamma_{ij}^k+\frac{1}{2F}\left(
F_i\d_j^k+F_j\d_i^k-F^k\rg_{ij}
\right)
\label{GammaHat}
\end{equation}
where we have used the notation $F_i=\partial_iF$ and $F^k=\rg^{jk}F_j$. 

Let us start with proof for first-order symbols.
With the
help of
$$
\widehat\nabla_i\phi=\nabla_i\phi-\frac{n\l}{2}\frac{F_i}{F}\,\phi
\qquad
\hbox{and}
\qquad
\widehat\nabla_iP_1^i=\nabla_iP_1^i+\frac{n(1-\d)}{2}\frac{F_iP_1^i}{F}
$$
and, using (\ref{nablaphi}), for every $P\in\cS_\d^1$ we find  
$$
\widehat\cQ_{\l,\m}(P)
=
\cQ_{\l,\m}(P)+\frac{n}{2}\left(\a(1-\d)-\l\right)\frac{F_iP^i}{F}.
$$
The equality $\widehat\cQ_{\l,\m}(P)=\cQ_{\l,\m}(P)$ is now equivalent to
(\ref{alpha}).

As for the second-order symbols, $P\in\cS_{2,\d}$, the proof involves the
calculation of $\widehat\nabla_i\widehat\nabla_j\phi$ and
$\widehat\nabla_i\widehat\nabla_jP^{k\ell}$, which is straightforward; it also
relies on the well-known transformation law \cite{Bes} of the Ricci tensor under
a conformal rescaling, $\widehat\rg=F\rg$, namely
\begin{equation}
\widehat\Ric
=
\Ric
-\frac{(n-2)}{2}\left(\frac{\nabla{dF}}{F}
-
\frac{3}{2}\,\frac{dF\otimes{}dF}{F^2}\right)
-
\half\left(\frac{\Delta{F}}{F}
-
\frac{(n-4)}{2}\,\frac{\Vert{}dF\Vert^2}{F^2}\right)\rg
\label{RicciHat}
\end{equation}
where $\Delta{F}=\rg^{ij}\nabla_i\partial_jF$ and
$\Vert{}dF\Vert^2=\rg^{ij}\partial_iF\partial_jF$. The scalar curvature
transforms accordingly as
\begin{equation}
\widehat{R}
=
\frac{R}{F} -(n-1)\left(
\frac{\Delta{F}}{F^2}+\frac{(n-6)}{4}\,\frac{\Vert{}dF\Vert^2}{F^3}
\right).
\label{RHat}
\end{equation}
Using  the formula (\ref{Q2}) as an Ansatz with undetermined
coefficients $\b_1,\ldots,\b_6$, a tedious calculation then shows that the
condition $\widehat\cQ_{\l,\m}(P)=\cQ_{\l,\m}(P)$ is equivalent to an
overdetermined linear system of $9$ equations for these coefficients. For
generic values of $\d$, the solution turns out to be unique and given by
(\ref{beta14}) and (\ref{beta56}).
\end{proof}

\goodbreak

%%%%%%%%%%%%%%%%%%%%%%%%%%%%%%%%%%%%%%%%%%%%%%%%%%%%%%%%%%%%%%%%%%%%%%%%%%%%%%%%%%%%
%%%%%%%%%%%%%%%%%%%%%%%%%%%%%%%%%%%%%%%%%%%%%%%%%%%%%%%%%%%%%%%%%%%%%%%%%%%%%%%%%%%%
\section{Applications}
\label{ApplicationSection}
%%%%%%%%%%%%%%%%%%%%%%%%%%%%%%%%%%%%%%%%%%%%%%%%%%%%%%%%%%%%%%%%%%%%%%%%%%%%%%%%%%%%
%%%%%%%%%%%%%%%%%%%%%%%%%%%%%%%%%%%%%%%%%%%%%%%%%%%%%%%%%%%%%%%%%%%%%%%%%%%%%%%%%%%%

We apply these results to the quantization of the geodesic flow on a conformally
flat manifold $(M,\rg)$, where, locally, $\rg_{ij}=F\,g_{ij}$ for
some smooth strictly positive function~$F$, i.e. to the quantization of the
quadratic polynomial $H=\rg^{ij}\xi_i\xi_j$ on~$T^*M$. We will furthermore
quantize the Hamiltonian $H=\rg^{ij}(\xi_i-A_i)(\xi_j-A_j)$, where
$A=A_i\,dx^i$ is a $U(1)$-connection, describing the motion of a charged
particle on a conformally flat manifold, minimally coupled to an electro-magnetic
field. We will also pay special attention to the resonant cases corresponding to the
table (\ref{TheArray}).

%%%%%%%%%%%%%%%%%%%%%%%%%%%%%%%%%%%%%%%%%%%%%%%%%%%%%%%%%%%%%%%%%%%%%%%%%%%%%%%%%%%%
\subsection{Quantization map}
%%%%%%%%%%%%%%%%%%%%%%%%%%%%%%%%%%%%%%%%%%%%%%%%%%%%%%%%%%%%%%%%%%%%%%%%%%%%%%%%%%%%

To define a quantization map out of the isomorphism $\cQ_{\l,\m}$, one
introduces (as usual) a real parameter $\hbar$ and replace the momenta $\xi_j$ by
their quantum substitutes~$i\hbar\xi_j$.  More specifically, let us consider a new
operator on symbols
$\cI_\hbar:\cS_\d\to(\cS_\d)^\bbC$ by
\begin{equation}
\cI_\hbar(P)(\xi)=P(i\hbar\,\xi).
\label{Ihbar}
\end{equation}
%We then propose the following
\begin{defi}
We will call conformally equivariant quantization the $\so(p+1,q+1)$-equivariant map
$\cQ_{\l,\m;\hbar}:\cS^2_\d\to(\cD_{\l,\m}^2)^\bbC$ defined by
\begin{equation}
\cQ_{\l,\m;\hbar}=\cQ_{\l,\m}\circ\cI_\hbar
\label{quantumMap}
\end{equation}
where $\cQ_{\l,\m}$ is as in (\ref{Q1},\ref{Q2}) and
$\cI_\hbar$ given by (\ref{Ihbar}).
\end{defi}

Let us recall that if 
\begin{equation}
\l+\m=1
\label{symmCond}
\end{equation}
there exists, for compactly-supported
densities, a $\Vect(M)$-invariant pairing $(\cF_\l)^\bbC\otimes(\cF_\m)^\bbC\to\bbC$
defined by
\begin{equation}
\varphi\otimes\psi\mapsto\int_M{\!\overline{\varphi}\,\psi}.
\label{pairing}
\end{equation}
The quantization map enjoys the following crucial property
\begin{pro}
The differential operators $\cQ_{\l,\m;\hbar}(P)$ defined by (\ref{quantumMap})
are sym\-metric (i.e., formally self-adjoint) for the pairing (\ref{pairing})
provided (\ref{symmCond}) holds.
\end{pro}
\begin{proof}
A more general version of this proposition has been proved in \cite{DLO}. In
our case, it can be proved directly. Easy calculation already gives, in any
coordinate system, the (formal) adjoints
$(P^{jk}\partial_j\partial_k)^*=\partial_j\partial_k\circ{}P^{jk}$, and
$(P^j\partial_j)^*=-\partial_j\circ{}P^j$. Using the expression
(\ref{Symbol}) in an adapted coordinate system, we find that
$\cQ_{\l,\m;\hbar}(P)$ is symmetric for any $P\in\cS^2_\d$ if and only if
$\a=\half$, $\b_1=1$ and $\b_2=0$. Returning to the values (\ref{alpha}) and
(\ref{beta14}) of the numeric coefficients, these conditions are satisfied if
$\l+\m=1$.
\end{proof}

%\goodbreak

%%%%%%%%%%%%%%%%%%%%%%%%%%%%%%%%%%%%%%%%%%%%%%%%%%%%%%%%%%%%%%%%%%%%%%%%%%%%%%%%%%%%
\subsection{Conformally equivariant Laplacian in the generic case}
%%%%%%%%%%%%%%%%%%%%%%%%%%%%%%%%%%%%%%%%%%%%%%%%%%%%%%%%%%%%%%%%%%%%%%%%%%%%%%%%%%%%

Consider the quadratic polynomial $H\in\cS_{2,\d}$ given in local coordinates on
$T^*M$ by
\begin{equation}
H=\rg^{ij}\xi_i\xi_j
\label{Hdelta}
\end{equation}
where $\rg=\rg_{ij}\,dx^i\otimes{}dx^j$ is a pseudo-Riemannian metric of
signature
$p-q$ on~$M$. 

\begin{pro}
\label{LapThm}
In the case $n\geq2$, and for $\l,\m$ fulfilling the condition (\ref{Reson}), the
quantization map (\ref{quantumMap}) yields the following expression:
\begin{equation}
\cQ_{\l,\m;\hbar}(H)=-\hbar^2\left(\Delta+C_{\l,\m}\,R\right)
\label{Geod}
\end{equation}
with
\begin{equation}
C_{\l,\m}=\frac{n^2\l(\m-1)}{(n-1)(n+2-2n\d)}
\label{C}
\end{equation}
where $\Delta$ is the Laplace operator and $R$ the scalar curvature
of~$(M,\rg)$.
\end{pro}
\begin{proof}
If $n\geq3$,
let us substitute the symbol $H$ given by (\ref{Hdelta}) into the formula
(\ref{Q2}). The second-order term of $\cQ_{\l,\m}(H)$ is nothing but the Laplace
operator $\Delta$. Since all covariant derivatives $\nabla_i\rg^{jk}$ vanish, we
are left with the scalar term $(\beta_5+n\beta_6)R$. The result follows from
(\ref{beta56}); note that the coefficient $\hbar^2$ comes from (\ref{quantumMap})
applied to the quadratic-homogeneous polynomial $H$.

In the case $n=2$, the result follows from (\ref{Q2S}) and (\ref{TraceSchwarz}).
\end{proof}

%%%%%%%%%%%%%%%%%%%%%%%%%%%%%%%%%%%%%%%%%%%%%%%%%%%%%%%%%%%%%%%%%%%%%%%%%%%%%%%%%%%%
\subsection{The Quantum Hamiltonian}
%%%%%%%%%%%%%%%%%%%%%%%%%%%%%%%%%%%%%%%%%%%%%%%%%%%%%%%%%%%%%%%%%%%%%%%%%%%%%%%%%%%%

In the special instance where $H\in\cS_0\cong\Pol(T^*M)$, the
Hamiltonian flow of $H$ projects onto the geodesics of $(M,\rg)$. 
Furthermore, in the most interesting case
\begin{equation}
\l=\m=\half
\label{half}
\end{equation}
naturally associated with geometric quantization, the operator
(\ref{Geod}) takes the form
\begin{equation}
\cQ_{\half,\half;\hbar}(H)
=
-\hbar^2\left(
\Delta-\frac{n^2}{4(n-1)(n+2)}\,R
\right).
\label{newfact}
\end{equation}
The self-adjoint operator (\ref{newfact}) on the
Hilbert space $\overline{\cF_\half}$ (the completion of the compactly
supported half-densities) is a natural new candidate for the quantized Hamiltonian
of the geodesic flow on a (pseudo-)Riemannian mani\-fold. None of the expressions
obtained in the literature by different methods of quantization (see, e.g.,
\cite{DO} for relevant references) corresponds to this one; all these expressions
therefore lack the conformal equivariance property (in the conformally flat case).

%%%%%%%%%%%%%%%%%%%%%%%%%%%%%%%%%%%%%%%%%%%%%%%%%%%%%%%%%%%%%%%%%%%%%%%%%%%%%%%%%%%%
\subsection{Minimal coupling and quantization}
\label{MinimaCcoupling}
%%%%%%%%%%%%%%%%%%%%%%%%%%%%%%%%%%%%%%%%%%%%%%%%%%%%%%%%%%%%%%%%%%%%%%%%%%%%%%%%%%%%

One can, as well, incorporate into the Hamiltonian (\ref{Hdelta}) additional terms
needed to describe electro-magnetic interaction. This is usually performed via the
so-called ``minimal coupling'' prescription to a
$U(1)$-connection, locally given by $A=A_i\,dx^i$. This procedure leads to a
Hamiltonian
$H\in\cS_\d$ of the form
\begin{equation}
H=\rg^{jk}(\xi_j-A_j)(\xi_k-A_k)
\label{minimalCoupling}
\end{equation}
on any pseudo-Riemannian manifold $(M,\rg)$. 

\begin{pro}
\label{MinimalPro}
In the case $n\geq2$, and for $\l,\m$ as in~(\ref{Reson}), the
quantization map~(\ref{quantumMap}) yields
\begin{equation}
\begin{array}{rcl}
\cQ_{\l,\m;\hbar}(H) &=&
\displaystyle
-\hbar^2
\rg^{jk}\Big(\nabla_j+\frac{i}{\hbar}A_j\Big)
\Big(\nabla_k+\frac{i}{\hbar}A_k\Big)
-\hbar^2C_{\l,\m}R\\[10pt]
&&+
i\hbar\,\displaystyle\frac{(1-\l-\m)}{(1-\d)}\,\rg^{jk}\nabla_jA_k
\end{array}
\label{QuantCouplingEq}
\end{equation}
where $C_{\l,\m}$ is given by (\ref{C}).
\end{pro}
The proof of the above proposition is completely analogous to that of Proposition
\ref{LapThm} and will be omitted.

Notice that the first line in (\ref{QuantCouplingEq}) corresponds to what is
called quantum minimal coupling in the physics literature. Thus, our conformally
equivariant quantization $\cQ_{\l,\m;\hbar}$ intertwines minimal coupling if and
only if condition (\ref{symmCond}) holds.

%%%%%%%%%%%%%%%%%%%%%%%%%%%%%%%%%%%%%%%%%%%%%%%%%%%%%%%%%%%%%%%%%%%%%%%%%%%%%%%%%%%%
\subsection{The resonant cases: the Yamabe operator}\label{ResonantCases}
%%%%%%%%%%%%%%%%%%%%%%%%%%%%%%%%%%%%%%%%%%%%%%%%%%%%%%%%%%%%%%%%%%%%%%%%%%%%%%%%%%%%

According to Theorem \ref{confEquivIsomThRes},
there exist pairs $(\l,\m)$ for which the modules~$\cD^2_{\l,\m}$
and~$\cS^2_{\m-\l}$ are isomorphic. However, we mentioned that this isomorphism
is not unique. But, imposing the condition (\ref{half}) for the module
$\cD^2_{\l,\m}$ enables us to look for the operators
$\cQ_{\l,\m}(H)$ which are symmetric (formally self-adjoint). 

\begin{pro}
\label{YamabeAndSons}
In each of the following resonant cases, there exists a unique isomorphism
$\cQ_{\l,\m;\hbar}$ for which the operator $\cQ_{\l,\m;\hbar}(H)$ is
symmetric:
\begin{eqnarray}
\cQ_{\frac{n-2}{2n},\frac{n+2}{2n};\hbar}(H) & = &
-\hbar^2\Big(\Delta-\frac{n-2}{4(n-1)}\,R\Big),
\label{YamabeOp}\\[5pt]
\cQ_{0,1;\hbar}(H) & = &
-\hbar^2\Delta,
\label{LaplaceOp}\\[5pt]
\cQ_{-\frac{1}{n},\frac{n+1}{n};\hbar}(H) & = &
-\hbar^2\Big(\Delta+\frac{1}{(n-1)(n+2)}\,R\Big).
\label{NewOp}
\end{eqnarray}
\end{pro}
The proof of the preceding proposition will be
given in Section \ref{ProofYamabeAndSons}.

We notice that the constraint (\ref{half}) selects only three (out of five)
resonances in~(\ref{TheArray}).

We recognize in (\ref{YamabeOp}) the so-called ``Yamabe'' operator and in 
(\ref{LaplaceOp}) the ordinary Laplace operator on functions.
At last, the operator (\ref{NewOp}) is a new $\so(p+1,q+1)$-equivariant Laplacian
which should be put quite on the same footing as the other two.

\medskip
\noindent
\textbf{Remarks:} 
(a) 
It is well known that the Yamabe operator (\ref{YamabeOp}) is the unique Laplace
operator which is invariant under conformal changes of metrics: $\rg\to F\,\rg$.
In this framework, the symbol $H\in\cS_{2,\d}$ given by (\ref{Hdelta})
with $\d=\frac{2}{n}$ is also invariant under conformal changes of metrics.

(b)
In contradistinction with the operator (\ref{newfact}), the conformal
Laplacians (\ref{YamabeOp},\allowbreak \ref{LaplaceOp},\ref{NewOp}) cannot serve as
self-adjoint quantum-mechanical operators on a Hilbert space since
$\l\not=\m$. 

(c) It is worth mentioning that the numerical coefficients in front of the scalar
curvature in (\ref{YamabeOp},\ref{LaplaceOp},\ref{NewOp}) actually correspond to
the expression (\ref{C}) that holds in the generic case.

%%%%%%%%%%%%%%%%%%%%%%%%%%%%%%%%%%%%%%%%%%%%%%%%%%%%%%%%%%%%%%%%%%%%%%%%%%%%%%%%%%%%
\subsection{The Sturm-Liouville operator}
\label{SturmLiouvilleSect}
%%%%%%%%%%%%%%%%%%%%%%%%%%%%%%%%%%%%%%%%%%%%%%%%%%%%%%%%%%%%%%%%%%%%%%%%%%%%%%%%%%%%

The operator $\cQ_{\l,\m}(H)$ for the symbol (\ref{Hdelta}) can be computed in
the case $n=1$ by (\ref{Q2One}); it appears to be still defined in the
resonant case, $\d=2$. 
In general, it does not yield an $\Sl(2,\bbR)$-equivariant
quantization map $\cQ_{\l,\l+2;\hbar}$ unless $\l=-\half$ and $\m=\frac{3}{2}$ (the
``Yamabe'' weights in (\ref{YamabeOp})). In this case one obtains a
special instance of Sturm-Liouville operator
\begin{equation}
\cQ_{-\half,\frac{3}{2};\hbar}(H)
=
-\hbar^2\Big(
\Delta -
\frac{S(\varphi)}{2\rg}
\Big)
\label{SturmLiouville}
\end{equation}
which can be interpreted as the Yamabe operator in the one-dimensional case.

\goodbreak

We now start the more technical part of our work in which we will provide the
proofs of the main theorems. We will derive the formul{\ae} for the
conformally equivariant isomorphism (\ref{confEquivIsom}) by means of the
algebraic theory of invariants.

%%%%%%%%%%%%%%%%%%%%%%%%%%%%%%%%%%%%%%%%%%%%%%%%%%%%%%%%%%%%%%%%%%%%%%%%%%%%%%%%%%%%
%%%%%%%%%%%%%%%%%%%%%%%%%%%%%%%%%%%%%%%%%%%%%%%%%%%%%%%%%%%%%%%%%%%%%%%%%%%%%%%%%%%%
\section{Euclidean invariant theory}
\label{WeylBrauerSection}
%%%%%%%%%%%%%%%%%%%%%%%%%%%%%%%%%%%%%%%%%%%%%%%%%%%%%%%%%%%%%%%%%%%%%%%%%%%%%%%%%%%%
%%%%%%%%%%%%%%%%%%%%%%%%%%%%%%%%%%%%%%%%%%%%%%%%%%%%%%%%%%%%%%%%%%%%%%%%%%%%%%%%%%%%

In this section we will introduce a Lie algebra of differential operators acting
on the space of symbols $\cS_{\m-\l}$ and commuting with the canonical action of
the Euclidean algebra. The associated universal enveloping algebra will provide
us with the ingredients needed to construct the conformally equivariant map
$\cQ_{\l,\m}$ (see \cite{DLO} for the abstract theory of conformally equivariant
quantization) on a conformally flat
$n$-dimensional manifold. Throughout this section we will assume $n\geq2$.

\goodbreak

%%%%%%%%%%%%%%%%%%%%%%%%%%%%%%%%%%%%%%%%%%%%%%%%%%%%%%%%%%%%%%%%%%%%%%%%%%%%%%%%%%%%
\subsection{The Weyl-Brauer Theorem}
%%%%%%%%%%%%%%%%%%%%%%%%%%%%%%%%%%%%%%%%%%%%%%%%%%%%%%%%%%%%%%%%%%%%%%%%%%%%%%%%%%%%

Consider first the space of polynomials $\bbC[\xi_1,\ldots,\xi_n]$ with the
canonical action of the orthogonal Lie algebra $\so(p,q)$ with $p+q=n$, 
generated by
%the vector fields 
$
X_{ij}=\xi_i{\partial}/{\partial \xi^j}-\xi_j{\partial}/{\partial \xi^i}
$ 
(cf.~(\ref{ConfAlg})). A classical theorem \cite{Wey,Bra}
states that the commutant $\so(p,q)^!$ in the space $\End(\bbC[\xi_1,\ldots,\xi_n])$
is the  associative algebra generated by:
\begin{equation}
\rR=\xi^i\xi_i,\qquad
\rE=\xi_i\frac{\partial}{\partial \xi_i}+\frac{n}{2},\qquad
\rT=\frac{\partial}{\partial \xi^i}\frac{\partial}{\partial \xi_i}
\label{sl2}
\end{equation}
whose commutation relations are those of $\Sl(2,\bbR)$.
%\begin{equation}
%[\rR,\rE]=-2\rR,\qquad
%[\rE,\rT]=-2\rT,\qquad
%[\rT,\rR]=4\rE.
%\label{sl2.Comm}
%\end{equation}
We will find it useful to deal with the Euler operator
\begin{equation}
\Euler=\rE-\frac{n}{2}.
\label{Euler}
\end{equation}
This algebra is, in fact, the universal enveloping algebra  $U(\Sl(2,\bbR))$.
%If $n=1$, it is isomorphic to the quotient algebra
%$U(\Sl(2,\bbR))/\langle\cC+\frac{3}{4}\rangle$ where
%\allowbreak
%$\cC=\rE^2-\half(\rR\,\rT+\rT\,\rR)$ is the Casimir operator of
%$\Sl(2,\bbR)$.

\medskip
\noindent
\textbf{Remarks:}

(a) Straightforward computation yields the explicit formul{\ae}:
\begin{equation}
\matrix{
\rR(P_k)^{i_1\ldots i_kij}
\hfill
&=&
P_k^{(i_1\ldots i_k}g^{ij)},\hfill\cr\noalign{\smallskip}
\Euler(P_k)^{i_1\ldots i_k}
\hfill
&=&
k\,P_k^{i_1\ldots i_k},\hfill\cr\noalign{\smallskip}
\rT(P_k)^{i_1\ldots i_{k-2}}
\hfill
&=&
k(k-1)g_{ij}P_k^{iji_1\ldots
i_{k-2}},\hfill\cr
}
\label{RET}
\end{equation}
where round brackets denote symmetrization.

(b) It is worth noticing that the converse property holds:
$\Sl(2,\bbR)^!=U(\so(p,q))$ showing that $\Sl(2,\bbR)$ and
$\so(p,q)$ form a dual pair of Lie algebras.

%%%%%%%%%%%%%%%%%%%%%%%%%%%%%%%%%%%%%%%%%%%%%%%%%%%%%%%%%%%%%%%%%%%%%%%%%%%%%%%%%%%%
\subsection{The Lie algebra of Euclidean invariants}
%%%%%%%%%%%%%%%%%%%%%%%%%%%%%%%%%%%%%%%%%%%%%%%%%%%%%%%%%%%%%%%%%%%%%%%%%%%%%%%%%%%%

Consider then the space
of polynomials $\bbC[x^1,\ldots,x^n,\xi_1,\ldots,\xi_n]$ with the canonical action
of the Euclidean Lie algebra $\se(p,q)=\so(p,q)\ltimes\bbR^n$ generated by the
canonical lifts to $T^*\bbR^n$ of the vector fields $X_{ij}$ and $X_i$ given by
(\ref{ConfAlg}). We are thus looking for the commutant $\se(p,q)^!$ in 
$\End(\bbC[x^1,\ldots,x^n,\xi_1,\ldots,\xi_n])$.
The following propositions extend the Weyl-Brauer theorem.

\begin{pro}
\label{CommEpq}
(i) The $\Sl(2,\bbR)$-module structure on $\bbC[x^1,\ldots,x^n,\xi_1,\ldots,\xi_n]$
extends to a module structure for the semi-direct product
$\Sl(2,\bbR)\ltimes\rh_1$, where $\rh_1$ is the Heisenberg Lie algebra generated
by:
\begin{equation}
\rG=\xi^i\frac{\partial}{\partial x^i},\qquad
\rD=\frac{\partial}{\partial \xi_i}\frac{\partial}{\partial x^i},\qquad
\rL=\frac{\partial}{\partial x^i}\frac{\partial}{\partial x_i}.
\label{rh1}
\end{equation}
(ii) The commutant $\se(p,q)^!$ is the associative algebra generated by the
operators given in (\ref{sl2}) and (\ref{rh1}).
\end{pro}
\begin{proof} Consider the commutant $\so(p,q)^!$ in the space
$\End(\bbC[x^1,\ldots,x^n,\xi_1,\ldots,\xi_n])$. As in the proof of the
Weyl-Brauer theorem we identify these endomorphisms with polynomials 
$\bbC[x^1,\ldots,x^n,p_1,\ldots,p_n,\xi_1,\ldots,\xi_n,y^1,\ldots,y^n]$, where the
$p_i$ and $y^i$ are in duality with $x^i$ and $\xi_i$ respectively. According to
\cite{Wey} the $\so(p,q)$-invariant polynomials are generated by the ten (scalar)
products: $x_ix^i,p_ix^i,\ldots,y_iy^i$. These second-order polynomials form a
Poisson algebra isomorphic to $\ssp(4,\bbR)$, therefore
$\so(p,q)^!$ is isomorphic to (some quotient of) $U(\ssp(4,\bbR))$.

The commutant $\se(p,q)^!$ is the subalgebra of $\so(p,q)^!$ which is
invariant under translations generated by $\partial/\partial x^i$. This subalgebra
is clearly generated by $\xi^i\xi_i$, $\xi_iy^i$, $y_iy^i$,
$\xi^ip_i$, $y^ip_i$, $p^ip_i$, in other words by the operators (\ref{sl2}) and
(\ref{rh1}).
\end{proof}
\medskip
\noindent
\textbf{Remarks:}

(a) Again, one easily finds:
\begin{equation}
\matrix{
\rG(P_k)^{i_1\ldots i_ki}\hfill
&=&
\partial_jP_k^{(i_1\ldots i_k}g^{i)j},\hfill\cr\noalign{\smallskip}
\rD(P_k)^{i_1\ldots i_{k-1}}
\hfill
&=&
k\,\partial_iP_k^{ii_1\ldots i_{k-1}},\hfill\cr\noalign{\smallskip}
\rL(P_k)^{i_1\ldots i_k}
\hfill
&=&
g^{ij}\partial_i\partial_jP_k^{i_1\ldots i_k}.
\hfill\cr
}
\label{GDL}
\end{equation}

(b) If $n\geq3$, one has $\so(p,q)^!=U(\ssp(4,\bbR))$ and
$\ssp(4,\bbR)^!=U(\so(p,q))$. This is also a well known instance of duality
between the orthogonal and symplectic algebras.

We furthermore prove the following

\goodbreak

\begin{thm}
\label{Usp3}
The commutant $\se(p,q)^!$ is isomorphic to %the quotient algebra
$U({\Sl(2,\bbR)\ltimes\rh_1})/\cI$ where the ideal $\cI$ is as follows:
\hfill\break
(i) if $n=2$, the ideal $\cI$ is generated by
\begin{equation}
Z = 
(\cC+{\textstyle\frac{3}{2}})\,\rL
+
{\textstyle\frac{1}{4}}\big(\rD\,[\rG,\cC] + [\rG,\cC]\,\rD
-
\rG\,[\rD,\cC] - [\rD,\cC]\,\rG
\big),
\label{Rel2}
\end{equation}
where $\cC=\rE^2-\frac{1}{2}(\rR\,\rT+\rT\,\rR)$ is the Casimir of $\Sl(2,\bbR)$,
\hfill\break
(ii) if $n\geq3$, one has 
\begin{equation}
\cI = \{0\}.
\label{Rel3}
\end{equation}
\end{thm}
\begin{proof}
Again, we identify the generators (\ref{sl2},\ref{rh1}) with the six
quadratic polynomials given in the preceding proof. 

If $n\geq3$, one finds that these polynomials are functionally, hence algebraically
independent. Indeed, 
$d(\xi^i\xi_i)\wedge{}d(\xi_jy^j)\wedge\cdots\wedge{}d(p^kp_k)\neq0$. 

In the case $n=2$, any five distinct polynomials from the previous set of quadratic
polynomials turn out to be independent. One then checks that the operator given by $Z$ in
(\ref{Rel2}) vanishes identically. Moreover, $Z\in{}U({\Sl(2,\bbR)\ltimes\rh_1})$ is of
minimal degree (three). Working, as above, in terms of polynomials (principal symbols), one
immediately gets, by using the implicit functions theorem, that any other polynomial in
this ideal is a multiple of the symbol of $Z$.
\end{proof}
\medskip

We do not know whether the converse to Theorem \ref{Usp3} is true:
our conjecture is that $(\Sl(2,\bbR)\ltimes\rh_1)^!=U(\se(p,q))$ for $n\geq3$; in
other words is it true that $U(\se(p,q))^{!!}=U(\se(p,q))$~? Similar problems have
recently been investigated by A.A.~Kirillov \cite{Kir2}.

%%%%%%%%%%%%%%%%%%%%%%%%%%%%%%%%%%%%%%%%%%%%%%%%%%%%%%%%%%%%%%%%%%%%%%%%%%%%%%%%%%%%
%%%%%%%%%%%%%%%%%%%%%%%%%%%%%%%%%%%%%%%%%%%%%%%%%%%%%%%%%%%%%%%%%%%%%%%%%%%%%%%%%%%%
\section{Equation characterizing conformal equivariance}
\label{ConfEquivSection}
%%%%%%%%%%%%%%%%%%%%%%%%%%%%%%%%%%%%%%%%%%%%%%%%%%%%%%%%%%%%%%%%%%%%%%%%%%%%%%%%%%%%
%%%%%%%%%%%%%%%%%%%%%%%%%%%%%%%%%%%%%%%%%%%%%%%%%%%%%%%%%%%%%%%%%%%%%%%%%%%%%%%%%%%%

%%%%%%%%%%%%%%%%%%%%%%%%%%%%%%%%%%%%%%%%%%%%%%%%%%%%%%%%%%%%%%%%%%%%%%%%%%%%%%%%%%%%
\subsection{Equivariance with respect to the affine subalgebra}
%%%%%%%%%%%%%%%%%%%%%%%%%%%%%%%%%%%%%%%%%%%%%%%%%%%%%%%%%%%%%%%%%%%%%%%%%%%%%%%%%%%%

We first consider, for the sake of completeness, the case of the whole affine Lie
subalgebra of $\Vect(\bbR^n)$.

\begin{lem}
\label{AffAction}
The actions (\ref{Dinf}) and (\ref{VectActionSymb},\ref{LieDer}) of the affine Lie
algebra $\gl(n,\bbR)\ltimes\bbR^n$ on the modules
$\cD_{\l,\m}$ and $\cS_{\m-\l}$ for the local expressions
(\ref{DiffOp}) and (\ref{Symb}) coincide identically.
\end{lem}
\begin{proof}
The $\Vect(M)$-action (\ref{Dinf}) has the following form in local coordinates:
\begin{equation}
L_X^{\l,\m}(A)_\ell=
L^{\m-\l}_X(A_\ell) + (\mbox{higher order derivatives of } X)
\label{ExplicitAction}
\end{equation}
for $X\in\Vect(M)$. The affine Lie algebra being characterized by the property that
all second derivatives $\partial_i\partial_j X^k$ vanish,
(\ref{ExplicitAction}) implies that each coefficient of the operator~$A$ transforms
as a symbol of degree~$\ell$.
\end{proof}
From now on, we identify locally the operators and the symbols by using the
formula~(\ref{ExplicitAction}).

%%%%%%%%%%%%%%%%%%%%%%%%%%%%%%%%%%%%%%%%%%%%%%%%%%%%%%%%%%%%%%%%%%%%%%%%%%%%%%%%%%%%
\subsection{Action of the inversions on $\cD^k_{\l,\m}$}
%%%%%%%%%%%%%%%%%%%%%%%%%%%%%%%%%%%%%%%%%%%%%%%%%%%%%%%%%%%%%%%%%%%%%%%%%%%%%%%%%%%%

At this stage, we need an explicit formula for the action (\ref{ExplicitAction}) of
the inversions, generated by $\bar X_i$ (see (\ref{ConfAlg})), on the space of
differential operators.

In order to make calculations more systematic, let us introduce the following
useful notation
\begin{equation}
L_{\bar X}=\xi_i\otimes L_{\bar X_i}
\label{Notation}
\end{equation}
which captures all the structure of the Abelian subalgebra of inversions.
Experience proved that this operator is compatible with all algebraic structures
introduced so far.

\begin{lem}
The action of the inversions on $\cD^k_{\l,\m}$
takes, with the convention (\ref{Notation}), the following form:
\begin{equation}
L_{\bar X}^{\l,\m}(A)_\ell=
L^{\m-\l}_{\bar X}(A_\ell) + 
(\ell+1)\left(-\half\ell\,\rR\,\rT+2(\ell+n\l)\right)A_{\ell+1}
\label{VeryExplicitAction}
\end{equation}
for $\ell=0,1,\ldots,k$.
\label{VEA}
\end{lem}
\begin{proof} Standard calculation leads to the general expression:
$$
\matrix{
L_X^{\l,\m}(A)_\ell^{i_1\ldots i_\ell} &=& 
L^{\m-\l}_X(A_\ell)^{i_1\ldots i_\ell}\hfill\cr
\noalign{\smallskip}
&&\displaystyle 
-\frac{\ell+1}{2}
\sum_{s=1}^\ell{A_{\ell+1}^{iji_1\ldots\widehat{i_s}\ldots{}i_\ell}
\partial_i\partial_jX^{i_s}}
-(\ell+1)\l A_{\ell+1}^{ii_1\ldots i_\ell}\partial_i\partial_jX^j\hfill\cr
\noalign{\smallskip}
&&
+\;(\mbox{higher order derivatives of } X)\hfill
}
$$
for any $X\in\Vect(M)$. In the case of inversions, namely, if $X=\bar X_r$, one has:
\begin{equation}
\partial_i\partial_j\bar X_r^s=2\left(g_{ij}\delta_r^s-
\delta_i^s g_{jr}-\delta_j^s g_{ir}\right),
\label{Partial}
\end{equation}
where $g_{ij}$ are the components of the flat metric on $\bbR^n$ given in Section
\ref{ConfFlatStruct}. The previous formula, therefore, becomes:
$$
\matrix{
L_{\bar X_r}^{\l,\m}(A)_\ell^{i_1\ldots i_\ell} &=& 
L^{\m-\l}_{\bar X_r}(A_\ell)^{i_1\ldots i_\ell}\hfill\cr
\noalign{\smallskip}
&&\displaystyle 
-(\ell+1)
\sum_{s=1}^\ell{g_{ij}\,A_{\ell+1}^{iji_1\ldots\widehat{i_s}\ldots{}i_\ell}
\delta^r_{i_s}}\hfill\cr
\noalign{\smallskip}
&&\displaystyle 
+2(\ell+1)(\ell+n\l) A_{\ell+1}^{ri_1\ldots i_\ell}
\hfill
}
$$
Then, using (\ref{RET}), one finds that the second term in the sum
$\xi_rL_{\bar X_r}^{\l,\m}(A)_\ell$ is equal to
$-\half\ell(\ell+1)\rR\,\rT(A_{\ell+1})$. The third term in the same expression is
plainly proportional to the identity.
\end{proof}

%%%%%%%%%%%%%%%%%%%%%%%%%%%%%%%%%%%%%%%%%%%%%%%%%%%%%%%%%%%%%%%%%%%%%%%%%%%%%%%%%%%%
\subsection{Equivariance equation}
%%%%%%%%%%%%%%%%%%%%%%%%%%%%%%%%%%%%%%%%%%%%%%%%%%%%%%%%%%%%%%%%%%%%%%%%%%%%%%%%%%%%

It is now possible to derive the main equation that guarantees the equivariance of
the symbol map and the quantization map with respect to the inversions.

\begin{pro}
A linear map $\cQ_{\l,\m}:\cS^k_{\m-\l}\to\cD^k_{\l,\m}$ intertwines
the action of the inversions if and only if the following equation holds:
\begin{equation}
[\cQ_{\l,\m},L_{\bar X}^{\m-\l}]
=
\left(
-{\textstyle\frac{1}{2}}\rR\,\rT(\Euler-1)+2\Euler+2(n\l-1)
\right)\Euler\circ\cQ_{\l,\m}.
\label{TheQEquation}
\end{equation} 
\label{EquivEq}
\end{pro}
\begin{proof} The equivariance condition writes: 
$\cQ_{\l,\m}\circ L_{\bar X}^{\m-\l}=L_{\bar X}^{\l,\m}\circ\,\cQ_{\l,\m}$.
Applying then equation (\ref{VeryExplicitAction}) to this condition readily yields
the result.
\end{proof}

%%%%%%%%%%%%%%%%%%%%%%%%%%%%%%%%%%%%%%%%%%%%%%%%%%%%%%%%%%%%%%%%%%%%%%%%%%%%%%%%%%%%
%%%%%%%%%%%%%%%%%%%%%%%%%%%%%%%%%%%%%%%%%%%%%%%%%%%%%%%%%%%%%%%%%%%%%%%%%%%%%%%%%%%%
\section{Proofs of the main results}\label{proofsMainResults}
%%%%%%%%%%%%%%%%%%%%%%%%%%%%%%%%%%%%%%%%%%%%%%%%%%%%%%%%%%%%%%%%%%%%%%%%%%%%%%%%%%%%
%%%%%%%%%%%%%%%%%%%%%%%%%%%%%%%%%%%%%%%%%%%%%%%%%%%%%%%%%%%%%%%%%%%%%%%%%%%%%%%%%%%%

%%%%%%%%%%%%%%%%%%%%%%%%%%%%%%%%%%%%%%%%%%%%%%%%%%%%%%%%%%%%%%%%%%%%%%%%%%%%%%%%%%%%
\subsection{Locality of the $\so(p+1,q+1)$-equivariant maps}
%%%%%%%%%%%%%%%%%%%%%%%%%%%%%%%%%%%%%%%%%%%%%%%%%%%%%%%%%%%%%%%%%%%%%%%%%%%%%%%%%%%%

It should be emphasized that the isomorphism (\ref{confEquivIsom}) is necessarily
given by a differential map, namely (\ref{Symbol}). This fact is already guaranteed
by the equivariance with respect to the subalgebra $\bbR\ltimes\bbR^n$ generated
by homotheties and translations (which is a common subalgebra of
$\so(p+1,q+1)$ and $\Sl(n+1,\bbR)$), i.e. by the
\begin{pro}\cite{LO}
If $k\geq\ell$, any $\bbR\ltimes\bbR^n$-equivariant map
$\cS^k_\delta\to\cS^\ell_\delta$ is local.
\label{LecomteOvsienko}
\end{pro} 
By Peetre's theorem \cite{Pee} such maps are locally given by differential
operators.

\goodbreak

%%%%%%%%%%%%%%%%%%%%%%%%%%%%%%%%%%%%%%%%%%%%%%%%%%%%%%%%%%%%%%%%%%%%%%%%%%%%%%%%%%%%
\subsection{The Ansatz}
%%%%%%%%%%%%%%%%%%%%%%%%%%%%%%%%%%%%%%%%%%%%%%%%%%%%%%%%%%%%%%%%%%%%%%%%%%%%%%%%%%%%

We will use our previous results on the universal enveloping algebra
$U({\Sl(2,\bbR)\ltimes\rh_1})$ to determine an adequate Ansatz for the
quantization map $\cQ_{\l,\m}:\cS^k_{\m-\l}\to\cD^k_{\l,\m}$, which turns out to
be more convenient in our framework. But, an identical general Ansatz would apply
just as well to the symbol map.

Proposition \ref{LecomteOvsienko}, together with the generalized Weyl-Brauer theorem
\ref{CommEpq}, leads to the general form for a $\se(p,q)$-equivariant quantization
map
$\cQ_{\l,\m}:\cS^k_{\m-\l}\to\cD^k_{\l,\m}$ given by differential operators
$
\cQ_{\l,\m}=C_{r,e,g,d,\ell,t}\,
\rR^r\,\rE^e\,\rG^g\,\rD^d\,\rL^\ell\,\rT^t,
$
where $C_{r,e,g,d,\ell,t}$ are constant coefficients.

Imposing, furthermore, the equivariance of $\cQ_{\l,\m}$ with
respect to homotheties generated by $X_0$ from (\ref{ConfAlg}), 
one readily finds
$t=r+g+\ell$ and
obtains that any $\so(p+1,q+1)$-equivariant map
$\cQ_{\l,\m}:\cS^k_{\m-\l}\to\cD^k_{\l,\m}$ is of the form
\begin{equation}
\cQ_{\l,\m}=C_{r,e,g,d,\ell}\,
\rR_0^r\,\Euler^e\,\rG_0^g\,\rD^d\,\rL_0^\ell,
\label{Ansatz}
\end{equation}
where we have put
\begin{equation}
\rR_0=\rR\,\rT,\qquad
\rG_0=\rG\,\rT,\qquad
\rL_0=\rL\,\rT.
\label{Zero}
\end{equation}

We will also impose the natural normalization condition which demands that the
principal symbol be preserved:
\begin{equation}
C_{r,e,0,0,0}
=
\left\{\matrix{
1&\mbox{if\ }(r,e)=(0,0)\hfill\cr
0&\mbox{otherwise.}\hfill
}
\right.
\label{Normalization}
\end{equation}

%%%%%%%%%%%%%%%%%%%%%%%%%%%%%%%%%%%%%%%%%%%%%%%%%%%%%%%%%%%%%%%%%%%%%%%%%%%%%%%%%%%%
\subsection{Solving the equivariance equation}
%%%%%%%%%%%%%%%%%%%%%%%%%%%%%%%%%%%%%%%%%%%%%%%%%%%%%%%%%%%%%%%%%%%%%%%%%%%%%%%%%%%%

In the case of second order differential operators, which is the one this article
is devoted to, our Ansatz (\ref{Ansatz}) implies that
$\se(p,q)$-equivariant maps:

(a) $\cS_\delta^k\to\cS_\delta^{k-1}$ are linear combinations of $\rD$ and
$\rG_0$ for $k=1,2$;

(b) $\cS_\delta^k\to\cS_\delta^{k-2}$ are linear combinations of $\rD^2$ and
$\rL_0$ for $k=2$ (note that in this special case the
other operators taken from (\ref{Ansatz}), namely $\rG_0^2$ and
$\rG_0\rD$ are expressible in terms of the latter).

Furthermore, the monomials in $\rR_0$
vanish because of the normalization condition~(\ref{Normalization}); the terms
$\rR_0\,\rD,\rR_0\,\rG_0,\ldots$ are identically zero for
$k\leq2$. 

\goodbreak

\begin{pro}\label{Solving}
There exists a unique quantization map 
\begin{equation}
\cQ_{\l,\m}=
\Id+\gamma_1\rG_0
+\gamma_2\rD
+\gamma_3\Euler\rD
+\gamma_4\rL_0
+\gamma_5\rD^2
\label{Ansatz2}
\end{equation}
satisfying the equivariance
equation (\ref{TheQEquation}) provided 
condition (\ref{Reson}) holds; it is given by:
\begin{equation}
\matrix{
\displaystyle
\gamma_1=
\frac{n(\l+\m-1)}{2(n\delta-2)(n(\delta-1)-2)},
\hfill\cr\noalign{\smallskip}
\displaystyle \gamma_2=
\frac{\l}{1-\delta},
\hfill\cr\noalign{\smallskip}
\displaystyle \gamma_3=
\frac{1-\l-\m}{(\delta-1)(n(\delta-1)-2)},\hfill
\cr\noalign{\smallskip}
\displaystyle \gamma_4=
\frac{n\l\Big(2+(4\l-1)n+(2\l^2-\l\m-\m^2+2\m-1)n^2\Big)}
{2(n(\delta-1)-1)(n(2\delta-1)-2)(n\delta-2)(n(\delta-1)-2)},
\hfill\cr\noalign{\smallskip}
\displaystyle
\gamma_5=
\frac{n\l(n\l+1)}{2(n(\delta-1)-1)(n(\delta-1)-2)}.
\hfill\cr
}
\label{TheSolution}
\end{equation}
\end{pro}
\begin{proof} Let us compute the left hand side of the equation
(\ref{TheQEquation}) where the quantization map given by our Ansatz (\ref{Ansatz2}).
We need the commutators of the differential operators entering (\ref{Ansatz2})
with the Lie derivative $L_{\bar X_i}^\delta$ with respect to the generators 
$\bar X_i$ given by (\ref{ConfAlg}). Using the notation (\ref{Notation}) we
first prove the
\begin{lem}
The following commutation relations hold:
\begin{equation}
\matrix{
[\rR_0,L_{\bar X}^\delta] \hfill &=&\;\;\,0,
\hfill\cr
\noalign{\medskip}
[\Euler,L_{\bar X}^\delta] \hfill &=&\;\;\,0,
\hfill\cr
\noalign{\medskip}
[\rG_0,L_{\bar X}^\delta] \hfill &=&\;\;\, 2\rR_0(\Euler-n\delta), \hfill\cr
\noalign{\medskip}
[\rD,L_{\bar X}^\delta] \hfill&=& -2\rR_0+4\Euler^2
-2(n(\delta-1)+2)\Euler, \hfill\cr
\noalign{\medskip}
[\rL_0,L_{\bar X}^\delta] \hfill &=& -4\rR_0\rD
+8\Euler\rG_0+2(n(1-2\delta)-2)\rG_0,\hfill\cr
\noalign{\medskip}
[\rD^2,L_{\bar X}^\delta] \hfill &=& -4\rR_0\rD-2\rG_0
+8\Euler^2\rD
+4(n(1-\delta)-1)\Euler\rD.\hfill\cr
}
\label{CommutationRelations}
\end{equation}
\label{Commutation}
\end{lem}
\begin{proof}
One finds, using 
(\ref{VectActionSymb},\ref{LieDer},\ref{ConfAlg}),
$[\rD,L_{\bar X_i}^\delta]
=
-2\xi_i\rT+4\cE\partial_{\xi_i}-2n(\delta-1)\partial_{\xi_i}$. Then, the final
expression for $[\rD,L_{\bar X}^\delta]$ follows from the definition of the
operators $\rR_0$ and $\Euler$ given by (\ref{sl2},\ref{Zero},\ref{Euler}). The
other commutators in (\ref{CommutationRelations}) are derived in the same fashion
with the help of the commutation relations of the operators (\ref{sl2})
and~(\ref{rh1}).
\end{proof}
\medskip

Using the commutation relations (\ref{CommutationRelations}), we find
%\begin{equation}
$$
\matrix{
%\displaystyle
[\cQ_{\l,\m},L_{\bar X}^{\m-\l}] &=&
\;\;\;2\gamma_1(\rR_0\cE-n\delta\rR_0)
\hfill\cr\noalign{\smallskip}
%\displaystyle 
&&+2\gamma_2(-\rR_0+2(\cE^2-\cE)-n(\delta-1)\cE)
\hfill\cr\noalign{\smallskip}
%\displaystyle
&&+2\gamma_3(-\rR_0(\cE-1)+(2-n(\delta-1))(\cE^2-\cE))
\hfill\cr\noalign{\smallskip}
%\displaystyle
 &&+2\gamma_4(-2\rR_0\rD+4\cE\rG_0+(n(1-2\delta)-2)\rG_0)
\hfill\cr\noalign{\smallskip}
%\displaystyle 
&&+2\gamma_5(-\rG_0-2\rR_0\rD+4\cE^2\rD+2(n(1-\delta)-1)\cE\rD)
\hfill\cr
}
%\label{TheLHS}
$$
%\end{equation}
while the right hand side of (\ref{TheQEquation}) is given by

{\flushleft{
\hspace{0.5cm}$
\begin{array}{rl}
(-\half\rR_0(\Euler-1)+2\Euler\!\!&\!\!+\;2(n\l-1))\Euler\circ\cQ_{\l,\m}=\\
\noalign{\smallskip}
&(-\half\rR_0(\Euler-1)+2(\Euler+n\l-1))\Euler\hfill\\
\noalign{\smallskip}
&+2(\Euler+n\l-1)\Euler(\gamma_1\rG_0+\gamma_2\rD+\gamma_3\Euler\rD)\hfill
%\label{TheRHS}
\end{array}
$
}}

\bigskip
%\vskip 0,3cm

\noindent
since the extra terms, namely
$
(\Euler-1)\Euler(\gamma_1\rG_0
+\gamma_2\rD
+\gamma_3\Euler\rD
+\gamma_4\rL_0
+\gamma_5\rD^2)
$
and
$
\Euler(\gamma_4\rL_0
+\gamma_5\rD^2)
$
obviously vanish on the space of second order symbols.

\goodbreak

Now, the equivariance condition (\ref{TheQEquation}) amounts to equating the two
previous expressions. Identifying the coefficients of $\rR_0, \rG_0, \rD$ and the
scalar terms (of order one and two), respectively, one gets the following system
of linear equations:
\begin{equation}
\left\{
\matrix{
(2-n\delta)\gamma_1-(\gamma_2+\gamma_3)=-\half, 
\hfill\cr\noalign{\smallskip}
(n(1-2\delta)+2)\gamma_4-\gamma_5=n\l\gamma_1, 
\hfill\cr\noalign{\smallskip}
2(n(1-\delta)+1)\gamma_5=n\l(\gamma_2+\gamma_3), 
\hfill\cr\noalign{\smallskip}
(1-\delta)\gamma_2=\l, 
\hfill\cr\noalign{\smallskip}
(2+n(1-\delta))(\gamma_2+\gamma_3)=n\l+1. \hfill\cr
}
\right.
\label{TheSystem}
\end{equation}
The solution of this system is unique and given by (\ref{TheSolution}).
\end{proof}

\medskip
\noindent
\textbf{Example:} Proposition (\ref{Solving}) yields, in particular, the following
half-density quantization map:
\begin{equation}
\cQ_{\half,\half}
=
\Id
+\half\rD
+\frac{n}{8(n+1)(n+2)}\rL_0
+\frac{n}{8(n+1)}\rD^2.
\label{QuantHalfDensity}
\end{equation}

%%%%%%%%%%%%%%%%%%%%%%%%%%%%%%%%%%%%%%%%%%%%%%%%%%%%%%%%%%%%%%%%%%%%%%%%%%%%%%%%%%%%
\subsection{Proof of Theorem \ref{confEquivIsomTh}}\label{ProofMainThm}
%%%%%%%%%%%%%%%%%%%%%%%%%%%%%%%%%%%%%%%%%%%%%%%%%%%%%%%%%%%%%%%%%%%%%%%%%%%%%%%%%%%%%%%%%%%%%%%%%%%%%%%%%%%%%%%%%%%%%%%%%%%%%%%%%%%%%%%%%%%%%%%%%%%%%%%%%%%%%%%%%%%%%%%%

The $\so(p+1,q+1)$-equivariant quantization map (\ref{Ansatz2}) precisely coincides
with the expression (\ref{Symbol}), since, taking into account the formul{\ae}
(\ref{RET},\ref{GDL}), one easily establishes the correspondence between the
coefficients (\ref{TheSolution}) and (\ref{alpha},\ref{beta14}).

We have thus proved the existence of an isomorphism (\ref{confEquivIsom}) provided
the coefficients (\ref{TheSolution}) are well-defined, i.e. condition
(\ref{Reson}) holds. This proves part (i) of Theorem~\ref{confEquivIsomTh}.

Then the formula (\ref{Ansatz}) and the normalization condition
(\ref{Normalization}) insure that, up to a multiplicative constant, every
$\so(p+1,q+1)$-equivariant quantization map (\ref{confEquivIsom}) is, indeed, of
the form (\ref{Ansatz2}). The uniqueness of the quantization map (part (ii) of
Theorem~\ref{confEquivIsomTh}) immediately follows  from Proposition~\ref{Solving}.
%\end{proof}
%%%%%%%%%%%%%%%%%%%%%%%%%%%%%%%%%%%%%%%%%%%%%%%%%%%%%%%%%%%%%%%%%%%%%%%%%%%%%%%%%%%%
\subsection{Proof of Theorem \ref{confEquivIsomThRes}}\label{ProofResonThm}
%%%%%%%%%%%%%%%%%%%%%%%%%%%%%%%%%%%%%%%%%%%%%%%%%%%%%%%%%%%%%%%%%%%%%%%%%%%%%%%%%%%%

The system (\ref{TheSystem}) determines all $\so(p+1,q+1)$-equivariant linear maps
from $\cS^2_{\m-\l}$ to $\cD^2_{\l\m}$. In the resonant cases, this
system has, in general, no solution. However, solving it for
$\gamma_1,\ldots,\gamma_5$ and $\l$ as an extra indeterminate, one immediately
obtains the values of $\l$ and~$\m$ displayed in (\ref{TheArray}). 
%\end{proof}
\medskip

In doing so, one finds that the coefficient $\gamma_3$ remains undetermined for the
third resonance, and $\gamma_4$ for the rest.

%%%%%%%%%%%%%%%%%%%%%%%%%%%%%%%%%%%%%%%%%%%%%%%%%%%%%%%%%%%%%%%%%%%%%%%%%%%%%%%%%%%%
\subsection{Proof of Proposition \ref{YamabeAndSons}}\label{ProofYamabeAndSons}
%%%%%%%%%%%%%%%%%%%%%%%%%%%%%%%%%%%%%%%%%%%%%%%%%%%%%%%%%%%%%%%%%%%%%%%%%%%%%%%%%%%%

Returning to the basic system (\ref{TheSystem}) in the presence of resonances, we
easily find that the free parameter $\gamma_3$ (resp. $\gamma_4$) is uniquely
determined, in each resonant case where $\l+\m=1$, if we require that the operators
$\cQ_{\l,\m;\hbar}(P)$ be symmetric for all~$P\in\cS^2_\d$. In such cases, the
explicit expressions (\ref{YamabeOp},\ref{LaplaceOp},\ref{NewOp}) are obtained in
the same manner as in the proof of Proposition \ref{LapThm}.
%\end{proof}

%%%%%%%%%%%%%%%%%%%%%%%%%%%%%%%%%%%%%%%%%%%%%%%%%%%%%%%%%%%%%%%%%%%%%%%%%%%%%%%%%%%%
\subsection{Proof of Proposition \ref{ProQ1}}\label{ProofProQ1}
%%%%%%%%%%%%%%%%%%%%%%%%%%%%%%%%%%%%%%%%%%%%%%%%%%%%%%%%%%%%%%%%%%%%%%%%%%%%%%%%%%%%

Let us consider a homogeneous first-order polynomial $P\in\cS_{1,\d}$. From the
expression~(\ref{nablaphi}) of the covariant derivative of a tensor density, one
gets
$$
\nabla_iP^i=\partial_iP^i + (1-\d)\Gamma_iP^i.
$$
We then deduce from the formula (\ref{Symbol}) in an adapted coordinate system that
$$
\begin{array}{rcl}
\cQ_{\l,\m}(P)
&=&\displaystyle
P^i\partial_i+\frac{\l}{1-\d}\Big(\nabla_iP^i-(1-\d)\Gamma_iP^i\Big)\\[10pt]
&=&\displaystyle
P^i\partial_i-\l\,\Gamma_iP^i+\frac{\l}{1-\d}\nabla_iP^i\\[10pt]
&=&\displaystyle
P^i\nabla_i+\a\,\nabla_iP^i
\end{array}
$$
thanks to (\ref{alpha}) and (\ref{nablaphi}). The proof for zero-order polynomials
is trivial.

%%%%%%%%%%%%%%%%%%%%%%%%%%%%%%%%%%%%%%%%%%%%%%%%%%%%%%%%%%%%%%%%%%%%%%%%%%%%%%%%%%%%
\subsection{Proof of Theorem \ref{ThmQ2Int}}\label{ProofThmQ2}
%%%%%%%%%%%%%%%%%%%%%%%%%%%%%%%%%%%%%%%%%%%%%%%%%%%%%%%%%%%%%%%%%%%%%%%%%%%%%%%%%%%%

Consider now a homogeneous second-order polynomial $P\in\cS_{2,\d}$
and let us, again, use the formula (\ref{Symbol}) in adapted coordinates.

The operator~$\cQ_{\l,\m}(P)$ can be rewritten
intrinsically. Indeed, the highest order term retains the following expression
$$
\begin{array}{rcl}
P^{ij}\partial_i\partial_j
&=&
P^{ij}\nabla_i\nabla_j
+
\left(P^{jk}\Gamma^i_{jk}
+
2\l{}P^{ij}\Gamma_j\right)\nabla_i\\[4pt]
&&
+P^{ij}\left(\l^2\Gamma_i\Gamma_j+\l\partial_i\Gamma_j\right)
\end{array}
$$
which can be deduced from (\ref{nablaphi}). Let us notice that, in order to obtain
such a seemingly standard expression, we actually need to differentiate
$\l$-densities and tensor fields with values in the space of $\l$-densities.
Neither the latter formula, nor the following ones are common in differential
geometry.

The two first-order terms in $\cQ_{\l,\m}(P)$ read
$$
\begin{array}{rcl}
(\partial_jP^{ij})\partial_i
&=&
(\nabla_jP^{ij})\nabla_i-\left(P^{jk}\Gamma^i_{jk}+(1-\d)P^{ij}\Gamma_j\right)
\nabla_i\\[4pt]
&&
+\l(\nabla_iP^{ij})\Gamma_j
-\l{}P^{ij}\left(\Gamma^k_{ij}\Gamma_k
+(1-\d)\Gamma_i\Gamma_j\right)
\end{array}
$$
and
$$
\begin{array}{rcl}
\rg^{ij}\rg_{k\ell}(\partial_jP^{k\ell})\partial_i
&=&
\rg^{ij}\rg_{k\ell}(\nabla_jP^{k\ell})\nabla_i
-\rg^{ij}\rg_{k\ell}\left(
2P^{m\ell}\Gamma^k_{jm}-\d{}P^{k\ell}\Gamma_j
\right)\nabla_i\\[4pt]
&&
+\rg^{ij}\rg_{k\ell}\left(
\nabla_jP^{k\ell}
-2P^{m\ell}\Gamma^k_{jm}+\d{}P^{k\ell}\Gamma_j
\right).
\end{array}
$$
At last, the two zero-order terms in $\cQ_{\l,\m}(P)$ are as follows
$$
\begin{array}{rcl}
\partial_i\partial_jP^{ij}
&=&
\nabla_i\nabla_jP^{ij}
-2(1-\d)(\nabla_iP^{ij})\Gamma_j
-(\nabla_iP^{jk})\Gamma^i_{jk}\\[4pt]
&&
-P^{ij}\left(\partial_k\Gamma^k_{ij}+(1-\d)\partial_j\Gamma_j
-2\Gamma_{ik}^\ell\Gamma_{j\ell}^k
-(1-2\d)\Gamma_{ij}^k\Gamma_k
-(1-\d)^2\Gamma_i\Gamma_j
\right)
\end{array}
$$
and
$$
\begin{array}{rcl}
\rg^{ij}\rg_{k\ell}\partial_i\partial_jP^{k\ell}
&=&
\rg^{ij}\rg_{k\ell}\Big(\nabla_i\nabla_jP^{k\ell}
-
4(\nabla_jP^{\ell{}m})\Gamma^k_{im}
+
2\d(\nabla_jP^{k\ell})\Gamma_i
\\[4pt]
&&-
2P^{\ell{}m}\partial_i\Gamma_{jm}^k
+
\d{}P^{k\ell}\partial_i\Gamma_j\\[4pt]
&&
+
2P^{\ell{}m}(\Gamma_{im}^r\Gamma_{jr}^k
-
2\d\Gamma_{im}^k\Gamma_j)
+
2P^{mr}\Gamma_{im}^k\Gamma_{jr}^\ell
+
\d^2P^{k\ell}\Gamma_i\Gamma_j
\Big).
\end{array}
$$

Now, to obtain the final formula for the conformally equivariant
map, let us collect the above terms within the expression~(\ref{Symbol}) where
the coefficients
$\b_1,\ldots,\b_4$ are considered undetermined. We also need to use the Christoffel
symbols of the conformally flat metric $\rg=Fg$, namely
\begin{equation}
\Gamma_{ij}^k
=
\frac{1}{2F}\left(
F_i\d_j^k+F_j\d_i^k-F^kg_{ij}
\right)
\label{GammaConfFlat}
\end{equation}
where $g$ is some flat metric and $F^k=g^{jk}F_j$ (see (\ref{GammaHat})).

The second-order term we get is plainly $P^{ij}\nabla_i\nabla_j$. Then, the
first-order term in~$\cQ_{\l,\m}(P)$ is just given by the second line of (\ref{Q2})
if we impose the following conditions
$$
\left(1-\b_1+\frac{n}{2}\left(2\l-\b_1(1-\d)\right)\right)\frac{P^{ij}F_j}{F}
=
0
$$
and
$$
\left(
\half(1-\b_1)+\b_2\Big(\frac{n\d}{2}-1\Big)\right)
g^{ij}g_{k\ell}\,\frac{P^{k\ell}F_j}{F^2}
=
0
$$
for the extra non-intrinsic terms; these conditions are satisfied if and only
if $\b_1$ and~$\b_2$ are as in (\ref{beta14}).
As for the zero-order terms, we, again, have to rule out two non-intrinsic
terms, viz
$$
\left(\l\b_1-2\b_3\Big(1-\d+\frac{1}{n}\Big)\right)(\nabla_iP^{ij})\Gamma_j
=
0
$$
and
$$
\left(
\l\b_2+\frac{\b_3}{n}+\b_4\Big(2\d-1-\frac{2}{n}\Big)
\right)
g^{ij}g_{k\ell}(\nabla_iP^{k\ell})\Gamma_j
=
0.
$$
These conditions determine $\b_3$ and $\b_4$ in accordance with (\ref{beta14}).

We finally check that the remaining zero-order terms in
$\cQ_{\l,\m}(P)$ are as follows
\begin{equation}
\begin{array}{rl}
\displaystyle\frac{n^2\l(1-\m)}{2(1+n(1-\d))}\!\!\!\!\!&\left[
\displaystyle\frac{P^{ij}F_{ij}}{F}
-\frac{3}{2}\frac{P^{ij}F_iF_j}{F^2}
+\frac{1}{2+n(1-2\d)}g^{ij}g_{k\ell}\left(
\frac{P^{k\ell}F_{ij}}{F}\right.\right.\\[14pt]
&
\left.\left.\displaystyle-\half(2+n(\d-1))\frac{P^{k\ell}{F_iF_j}}{F^2}
\right)
\right]
\end{array}
\label{zeroOrderTerms}
\end{equation}
where $F_{ij}=\partial_i\partial_jF$.
At this stage, some more ingredients are needed, namely the Ricci
tensor and the scalar curvature for the the conformally flat metric $\rg=Fg$
with Christoffel symbols~(\ref{GammaConfFlat}). The corresponding
expressions can be easily deduced from (\ref{RicciHat}) and~(\ref{RHat}).

One checks that, in the
case $n\geq3$, the expression (\ref{zeroOrderTerms}) 
organizes as the combination
$\b_5P^{ij}R_{ij}+\b_6P^{ij}\rg_{ij}R$ where $\b_5$ and $\b_6$ are rigidly fixed
and coincide with (\ref{beta56}).

%%%%%%%%%%%%%%%%%%%%%%%%%%%%%%%%%%%%%%%%%%%%%%%%%%%%%%%%%%%%%%%%%%%%%%%%%%%%%%%%%%%%
\subsection{Proof of Theorems \ref{ThmQ2n1} and \ref{ThmQ2n2}}\label{ProofThmQ2n12}
%%%%%%%%%%%%%%%%%%%%%%%%%%%%%%%%%%%%%%%%%%%%%%%%%%%%%%%%%%%%%%%%%%%%%%%%%%%%%%%%%%%%

In the lower dimensional cases, $n=1$ and $n=2$, the proofs are similar to that of
the higher dimensional case $n\geq3$ given in Section
\ref{ProofThmQ2}. All the computations giving~$\cQ_{\l,\m}(P)$ are exactly the same
as above right up to the formula (\ref{zeroOrderTerms}).

Now, in order to give an intrinsic interpretation of (\ref{zeroOrderTerms}), we
need to resort to the definition (\ref{multiDimSchwarzian}) of the Schwarzian
derivative in the case $n=2$. The final formula (\ref{Q2S}) for $\cQ_{\l,\m}(P)$
then readily follows from the expression (\ref{RHat}) of the scalar curvature for
the conformally flat metric $\rg=Fg$.

In the case $n=1$, the formula (\ref{Q2One}) obtained exactly in the same way as
(\ref{Q2S}).

%%%%%%%%%%%%%%%%%%%%%%%%%%%%%%%%%%%%%%%%%%%%%%%%%%%%%%%%%%%%%%%%%%%%%%%%%%%%%%%%%%%%
%%%%%%%%%%%%%%%%%%%%%%%%%%%%%%%%%%%%%%%%%%%%%%%%%%%%%%%%%%%%%%%%%%%%%%%%%%%%%%%%%%%%
\section{Conclusion and outlook}\label{ConclusionOutlook}
%%%%%%%%%%%%%%%%%%%%%%%%%%%%%%%%%%%%%%%%%%%%%%%%%%%%%%%%%%%%%%%%%%%%%%%%%%%%%%%%%%%%
%%%%%%%%%%%%%%%%%%%%%%%%%%%%%%%%%%%%%%%%%%%%%%%%%%%%%%%%%%%%%%%%%%%%%%%%%%%%%%%%%%%%

In this work, we have taken a first step towards a conformally invariant
quantization, i.e., depending only on the conformal class of a pseudo-Riemannian
metric. This program is now achieved for the case of second-order symbols
and differential operators. The general case still remains to be tackled, however
computations seem much more intricate.

Our original idea was to relate geometric quantization and deformation
quantization in a somewhat novel fashion, namely by using, from the start,
equivariance with respect to some structural symmetry group (e.g. the conformal
group). In the conformally flat case, it has been proved \cite{DLO} that there
exists, for any order, a canonical quantization map equivariant with respect to the
action of the conformal group. This conformally flat case is particular and we have
been able to extend the second-order quantization map to the case of an arbitrary
pseudo-Riemannian manifold.

As a by-product, we have obtained a new quantization of the geodesic
flow~(\ref{newfact}) on the Hilbert space of half-densities. We have also related the
Yamabe operator to other conformally equivariant Laplacians on resonant modules of
densities, and derived the quantum version of minimal coupling in the same
framework. 

We have also chosen to put aside the cohomological content of many aspects of the
problem. It should be stressed that Lie-algebra cohomology proved useful in
earlier work \cite{DO,LO,Gar,Lec} on the modules of differential operators. The
resonances appearing in (\ref{Reson}) should thus certainly hide non-trivial
$\so(p+1,q+1)$-cohomology classes.

Let us finish by mentioning a crucial property of the conformal algebra which was of
central importance in our work. The Lie algebra $\so(p+1,q+1)$ is a maximal Lie
subalgebra of $\Vect(\bbR^n)$ in the sense that any larger subalgebra is
infinite-dimensional (see \cite{BL}). This property implied the uniqueness of the
isomorphisms of the modules of differential operators and symbols under study.
Recall that the same is true for the projective Lie algebra
$\Sl(n+1,\bbR)$.

\bigskip

\textbf{Acknowledgments:} We are indebted to J.-L.~Brylinski and A.A.~Kirillov  for
enlightening discussions. We are also grateful to P.~Lecomte, Z.J.~Liu,
O.~Ogievetsky, C.~Roger and S.~Loubon-Djounga for helpful suggestions.

%\newpage

%%%%%%%%%%%%%%%%%%%%%%%%%%%%%%%%%%%%%%%%%%%%%%%%%%%%%%%%%%%%%%%%%%%%%%%%%%%%%%%%%%%%
%%%%%%%%%%%%%%%%%%%%%%%%%%%%%%%%%%%%%%%%%%%%%%%%%%%%%%%%%%%%%%%%%%%%%%%%%%%%%%%%%%%%

\end{document}